\documentclass[10pt,twoside]{article}
\usepackage{mathrsfs}
\usepackage{amssymb}
\usepackage{amsmath,amsbsy,amsfonts}
\usepackage[usenames]{color}
\usepackage[dvipsnames,svgnames,x11names]{xcolor}
\usepackage[english]{babel}
\usepackage[utf8]{inputenc}
\usepackage{authblk}
\usepackage{verbatim}
\newtheorem{lemma}{Lemma}[section]
\newtheorem{proposition}{Proposition}[section]
\newtheorem{theorem}{Theorem}[section]

\newtheorem{remark}{Remark}[section]

\newcommand{\I}{{\mathcal I}}

\newcommand{\R}{{\mathbb R}}
\newcommand{\N}{{\mathbb N}}

\renewcommand{\S}{{\mathcal S}}
\renewcommand{\phi}{{\varphi}}

\newcommand{\eps}{{\varepsilon}}
\renewcommand{\epsilon}{{\varepsilon}}

\numberwithin{equation}{section}

\let\Section=\section
\def\section{\setcounter{equation}{0}\Section}

\title{Multiplicity results for critical fractional Ambrosetti-Prodi type system with nonlinearities interacting with the spectrum}

\author[1]{Eduardo. H. Caqui}
\author[2]{Sandra M. de S. Lima}
\author[3]{F\'abio R. Pereira}
\affil[1]{Departamento de Ciencias Sede Bre$\widetilde{n}$a, Universidad Privada del Norte, Av. Tingo María 1122, Cercado de Lima, Lima, Peru.
\textit{eduardocaqui11@gmail.com}}
\affil[2]{Departamento de Ciências Exatas, Biológicas e da Terra, INFES-UFF, Santo Antônio de Pádua - RJ, 28470-000, Brazil. \textit{sandramsl@id.uff.br}}
\affil[3] {Departamento de Matem\'atica,
Instituto de Ci\^{e}ncias Exatas,
Universidade Federal de Juiz de Fora, 30161-970\\
Juiz de Fora - MG, Brazil.\textit{fabio.pereira@ufjf.edu.br}\thanks{F. R. Pereira was supported partially by FAPEMIG/Brazil (RED-00133-21) and FAPEMIG/Brazil (CEX APQ 04528/22).}}

\date{}

\begin{document}
\maketitle

\begin{abstract}
We investigated the existence of solutions for a class of Ambrosetti-Prodi type systems involving the fractional Laplacian
operator and with nonlinearities reaching critical growth and interacting, in some sense, with the spectrum of the operator. The resonant case in $\lambda_{k,s}$ for $k>1$ is also investigated.
\end{abstract}

{\scriptsize{\bf 2000 Mathematics Subject Classification:} 35B06, 35B09, 35J15, 35J20.}

{\scriptsize{\bf Keywords:} Ambrosetti-Prodi problem, Fractional systems, resonance, critical growth.}

	\section{Introduction}

Let $s\in (0,1)$, $N>2s$ and $\Omega\subset\R^N$ is a bounded smooth domain.
In this paper we study the possibility of existence of solutions for the following critical fractional system
	
\begin{equation}\label{eq:1.1}
\small{\left\{
\begin{aligned}
(-\Delta)^s u &=au+bv+\frac{\alpha }{\alpha +\beta}{u_+}^{\alpha-1}{v_+}^\beta + \xi_1 {u_+}^{\alpha+\beta-1} + f && \text{in $\Omega$,}\\
(-\Delta)^s v &=bu+cv+\frac{\beta }{\alpha +\beta}{u_+}^\alpha{v_+}^{\beta-1} + \xi_2 {v_+}^{\alpha+\beta-1}  + g && \text{in $\Omega$,}\\
u=&\, v=0 && \text{in $\R^N\setminus\Omega$,}
\end{aligned}
\right.}
\end{equation}
where \[ (-\Delta)^su(x):=C(N,s)\lim_{\eps\searrow 0}\int_{\R^N\setminus B_\eps(x)}\frac{u(x)-u(y)}{|x-y|^{N+2s}}\,dy, \quad\,\, x\in\R^N,
\] is the fractional laplacian operator with
$C(N,s) =\Big(\displaystyle\int_{\R^N} \frac{1- cos(\zeta_1)}{|\zeta|^{N+2s}} d\zeta \Big)^{-1}$ a  positive dimensional constant,
$\alpha,\beta>1$ are real
constants such that the sum $\alpha+\beta$ is the fractional critical Sobolev  exponent $ 2^*_s:=\frac{2N}{N-2s}$, $\xi_1, \xi_2 \geq 0$ and the forcing terms $f$ and  $g$ are of the form $f=t\phi_{1,s}+f_1$ and  $g=r\phi_{1,s}+g_1$, in such a way that the pair $(t,r) \in
\mathbb{R}^2$, $f_1,g_1\in L^q (\Omega)$ for some $q> \frac{N}{2s}$ and $\int_{\Omega} f_1\phi_{1,s} dx= \int_{\Omega} g_1\phi_{1,s}
dx=0$ with $\phi_{1,s}$ the positive eigenfunction associated with the first eigenvalue $\lambda_{1,s}$ of the operator $(-\Delta)^s$ with homogeneous Dirichlet boundary condition.

\vskip4pt
\noindent
With the above decomposition, in order to state and compare our
results to the scalar case, it is convenient to rewrite system
(\ref{eq:1.1}) as
\begin{equation}\label{c1.1}
\left\{
\begin{array}{lc}
(-\overrightarrow{\Delta})^s U =AU +  \nabla F(U) + T\phi_{1,s}+F_1
&\; \text{in $\Omega$,} \\
U=0 & \;\text{in $\R^N\setminus\Omega$,} \end{array}\right.
\end{equation}
where $U=\left( \begin{array}{c}    u \\
v \\
\end{array}
\right),
(-\overrightarrow{\Delta})^s U=\left(
\begin{array}{cc}
(-\Delta)^s u & 0 \\
0 & (-\Delta)^s v \\
\end{array}
\right), A=\left(
\begin{array}{cc}
a &  b \\
b & c \\
\end{array}
\right)
\in M_{2 \times 2}(\mathbb{R}),$ $\nabla$ is the gradient operator,
$F(U)=\displaystyle \frac {1}{\alpha +\beta } \left({u_+}^{\alpha}{v_+}^\beta + \xi_1 {u_+}^{\alpha+\beta}+ \xi_2 {v_+}^{\alpha+\beta}\right),$\\$T=\left( \begin{array}{c}
t \\
r \\
\end{array}
\right)
$ and $F_1=\left(
\begin{array}{c}
f_1 \\
g_1 \\
\end{array}
\right).$\\

Let $\mu_1,\mu_2$ be real eigenvalues of the symmetric matrix $A$, which will assume $\mu_1 \leq \mu_2.$ Thus, it is verified that
$\mu_1 |U|^2\leq (AU,U)_{\R^2}\leq \mu_2 |U|^2,\quad\text{for all $U:=(u,v) \in\R^2$.} $
The interaction of these eigenvalues with the spectrum of the $ (-\Delta)^s$ will play an important role in the study of existence of the solutions.

We recall that Ambrosetti and Prodi \cite{ap} in 1972, studied the following boundary value problem

\begin{equation}
\small{\left\{
\begin{aligned}\label{eq:0}
-\Delta u &=f(u)+g(x) && \text{in $\Omega$,}\\
u=&0 && \text{on $\partial\Omega$,}
\end{aligned}
\right.}
\end{equation}
where $g \in C^{0,\alpha}(\overline{\Omega})$ with $\alpha \in (0, 1), f \in C^2(\R)$ such that $f(0) = 0, f''
(t) > 0 $ for all $t \in \R$ and
$$0< \lim_{t\rightarrow -\infty} f'(t)< \lambda_1< \lim_{t\rightarrow +\infty} f'(t)<\lambda_2,$$ where $0<\lambda_1<\lambda_2 \leq  \ldots  \leq \lambda_k \ldots$ denote the eigenvalues of $(-\Delta, H_0^1(\Omega)).$
The authors showed that there exists in $C^{0,\alpha}(\overline{\Omega}),$ a closed connected $C^1$ manifold $M_1$ of codimension 1 which splits the space into two connected components $M_0$ and $M_2$  such that, if $g \in M_0$, the problem (\ref{eq:0}) has no solution; if $g \in M_1$, the problem (\ref{eq:0}) has exactly one solution and if  $g \in M_2$, the problem (\ref{eq:0}) has exactly two solutions.
After the pioneering work by Ambrosetti and Prodi \cite{ap}, many existence and multiplicity results have been investigated in different directions.
In particular, Ruf and Srikanth
(in \cite{rs}) established a multiplicity result for the local subcritical problem $ - \Delta  u = \lambda u + u_+^p + f ( x ) \;\mbox{in} \; \Omega, \; u = 0 \; \mbox{on}\;
\partial \Omega $ provided that the non-homogeneous term $f$ has the form $f (x) = h(x) + t\varphi_1 (x)$ ($h \in L^r (\Omega)$ with $r > N$), $\lambda$
is not an eigenvalue of $(- \Delta, H_0^1(\Omega))$ and $t > T,$ for some sufficiently large number $T = T (h).$
Still in the local scalar case, but with nonlinearity in the critical growth $(p=2^* -1)$, the problem above mentioned has been studied by De Figueiredo and Jianfu (in \cite{djj}), the authors proved the
existence of two solutions when $N > 6.$ This result was extended by Calanchi and Ruf \cite{cr} using the technique developed in \cite{gr}. Works related to this subject in the local scalar case, we recommend \cite{av} and in the nonlocal operators situation, \cite{AMBROSIO} and \cite{fx}
(see references therein).
For the  critical system in the local operators situation, problem (\ref{eq:1.1}) was studied,
for instance, in \cite{fd} and \cite{bru} when $\mu_2 < \lambda_1$ and by \cite{fp} in the uncoupled case.
For the fractional subcritical system, (\ref{eq:1.1}) was studied, for instance, in \cite{frp}.

The purpose of this work is to prove the existence of solutions for the class of nonlocal gradient systems of elliptic equations (\ref{eq:1.1}) involving critical nonlinearities on the hypothesis of an interaction of the eigenvalues $\mu_1, \mu_2$ of the matrix $A$ with eigenvalues of the fractional Laplacian operator $(-\Delta)^s.$
When $\mu_2 < \lambda_{1,s}$, this system belongs to the class of the so called Ambrosetti-Prodi type problems \cite{ap} which have been studied by several authors
in the last decades with different approaches.

Problem (\ref{eq:1.1}) is an extension to systems involving fractional Laplacian operator of the equation considered in \cite{rs}, \cite{djj} and \cite{cr}, in which (\ref{eq:1.1})
was studied in the local operators case ($s=1$) and nonlocal operators ($0<s<1$) in  \cite{AMBROSIO} (see \cite{fx} also) and with the particular matrix $ A=\left(
\begin{array}{cc}
\lambda &  0 \\
0 & \lambda \\
\end{array}
\right)
\in M_{2 \times 2}(\mathbb{R}).$
In this paper, we complement the results achieved in \cite{frp}, proving that the system (\ref{eq:1.1}) (or (\ref{c1.1})) has at least two solutions for sufficiently large values of parameters $(t,r)$, the first solution is negative and obtained explicitly depending on the non-homogeneous terms $f$ and $g$. The second solution is obtained via the Mountain Pass Theorem
when $\mu_2<\lambda_{1,s}$, or applying the Linking Theorem in the case $\lambda_{k,s}<\mu_1\leq\mu_2<\lambda_{k+1,s}$ if $k\geq1$.
The resonant case $\lambda_{k,s}=\mu_1$ for $k>1$ is also treated here.\\
Finally, we should point out that the corresponding local problem governed by the standard Laplacian operator can be recovered by letting $s\rightarrow 1.$

To show the existence of solution, difficulties arise when we consider fractional operators. As we know, in \cite{cr}, the  approximate eigenfunctions technique was used to facilitate the estimates of the energy functional associated with the local scalar problem in the space $H_0^1(\Omega)$ (for local critical systems, also see $\cite{bru}$).  However, as noted in $\cite{mmp}$, in the nonlocal case, it is not possible to employ any more the same idea as in \cite{cr} or $\cite{bru}$, since $u$ and $v$ are not orthogonal in the fractional space $X_0^s(\Omega)$ even though they have disjoint supports.
On the other hand, further complications arise due to the presence of the mathematical term
$F(u,v)= \displaystyle \frac {1}{\alpha +\beta } \left[{u_+}^{\alpha}{v_+}^\beta + \xi_1 {u_+}^{\alpha+\beta}+ \xi_2 {v_+}^{\alpha+\beta}\right]$ that
includes either an uncoupled or a coupled  nonlinearity.

Due to these obstacles, we develop similar techniques to these known for the Laplacian operator.

It is important to point out that, with the aid of \cite{MD}, our
results are still valid for the general case $\nabla F(u,v)$ when $F$ is a $(\alpha+\beta)-$homogeneous nonlinearity, which includes a
larger class of functions. \\

The proof of the Theorem below follows arguments as in \cite{frp}, so we will omit some details.

\begin{theorem}[Existence of a negative solution]\label{teo1}
Let $A \in M_{2 \times 2} (\R)$ be a symmetric matrix such that
\begin{eqnarray}
&&\det({\lambda_{j,s}}\; I-A) \neq 0, \; \forall \; j=1,2, \dots . \label{eq:1.3}
\end{eqnarray}
Assume that
$ F_1 = ( f_1 , g_1 ) \in L^q (\Omega) \times
L^q (\Omega)$ for some $q> \frac{N}{2s}$
and consider
\begin{equation}\nonumber
\textbf{\textit{R}}=\left\{ (t,r) \in \R^2  :   \begin{array}{c} \begin{aligned} br + (\lambda_{1,s} -c)t &< \eta \det({\lambda_{1,s}}I-A)\; \; \mbox{and}\\
(\lambda_{1,s}-a)r + bt &< \vartheta \det({\lambda_{1,s}}I-A)\} \end{aligned}  \end{array} \right\}.
\end{equation}
Then there exist $\eta, \vartheta \ll 0$ such that system (\ref{c1.1}) has a solution
$(u_{T} , v_{T})$ (with $u_{T}<0$ and $v_{T}<0$ in $\Omega$) for every
$T \in \textbf{\textit{R}}.$
\end{theorem}

\begin{remark}\label{obs0}
Suppose that $\; \det({\lambda_{1,s}}I-A) > 0$ and
\begin{eqnarray}
&& \lambda_{1,s} > \max \{a,c \}, \label{eq:1.4}
\end{eqnarray}
then the set $\textbf{\textit{R}}$ is a region between lines satisfying:\\
$(i)$ If $b=0$,
$$\textbf{\textit{R}}=
(-\infty \, ,\eta \, \frac{\lambda_{1,s}-c}{\det({\lambda_{1,s}}I-A)}) \times (-\infty \, , \vartheta \, \frac{\lambda_{1,s}-a}{\det({\lambda_{1,s}}I-A)}) \subset \R^2.$$
$(ii)$ If $b>0$,
\begin{equation}\nonumber
\textbf{\textit{R}}= {\scriptsize{ \left\{ (t,r) \in \R^2  :  \begin{array}{c} \begin{aligned} r  &< \eta \frac{\det({\lambda_{1,s}}I-A)}{b} - \frac{ (\lambda_{1,s} -c)}{b}t \; \; \mbox{and}\\
r  &< \vartheta \frac{\det({\lambda_{1,s}}I-A)}{\lambda_{1,s}-a} - \frac{b}{\lambda_{1,s}-a}t \end{aligned}   \end{array} \right\}}}.
\end{equation}
$(iii)$ If $b<0$,
\begin{equation}\nonumber
\textbf{\textit{R}}= {\scriptsize{ \left\{ (t,r) \in \R^2  :   \begin{array}{c} \begin{aligned} r  &> \eta \frac{\det({\lambda_{1,s}}I-A)}{b} - \frac{ (\lambda_{1,s} -c)}{b}t \; \; \mbox{and}\\
r  &< \vartheta \frac{\det({\lambda_{1,s}}I-A)}{\lambda_{1,s}-a} - \frac{b}{\lambda_{1,s}-a}t \end{aligned}  \end{array} \right\}}}.
\end{equation}

On the other hand, if  $\; \det({\lambda_{1,s}}I-A) > 0$ and
\begin{eqnarray}
&& \lambda_{1,s} < \min \{a,c \}, \label{eq:1.40}
\end{eqnarray}
then the set $\textbf{\textit{R}}$ satisfies:\\
$(i)$ If $b=0$,
$$\textbf{\textit{R}}= (\eta \, \frac{\lambda_{1,s}-c}{\det({\lambda_{1,s}}I-A)} \, , +\infty) \times (\vartheta \, \frac{\lambda_{1,s}-a}{\det({\lambda_{1,s}}I-A)}\, , +\infty) \subset \R^2.$$
$(ii)$ If $b>0$,
\begin{equation}\nonumber
\textbf{\textit{R}}= {\scriptsize{ \left\{ (t,r) \in \R^2  :  \begin{array}{c} \begin{aligned} r  &< \eta \frac{\det({\lambda_{1,s}}I-A)}{b} - \frac{ (\lambda_{1,s} -c)}{b}t \; \; \mbox{and}\\
r  &> \vartheta \frac{\det({\lambda_{1,s}}I-A)}{(\lambda_{1,s}-a)} - \frac{b}{(\lambda_{1,s}-a)}t \} \end{aligned}  \end{array} \right\}}}.
\end{equation}
$(iii)$ If $b<0$,
\begin{equation}\nonumber
\textbf{\textit{R}}= {\scriptsize{ \left\{ (t,r) \in \R^2  :  \begin{array}{c} \begin{aligned} r  &> \eta \frac{\det({\lambda_{1,s}}I-A)}{b} - \frac{ (\lambda_{1,s} -c)}{b}t \; \; \mbox{and}\\
r  &> \vartheta \frac{\det({\lambda_{1,s}}I-A)}{(\lambda_{1,s}-a)} - \frac{b}{(\lambda_{1,s}-a)}t \}  \end{aligned}  \end{array} \right\}}}.
\end{equation}
Note that, since $\det({\lambda_{1,s}}\; I-A) \neq 0$, the lines that define the region \textbf{\textit{R}} are not parallel.
Moreover, if $\det({\lambda_{1,s}}\; I-A) < 0$ a similar result can be obtained as in the Remark \ref{obs0}.
\end{remark}

The following are the main results of the paper.

\begin{theorem}\label{teo2}
Assume that  $N>6s$, $\xi_1  ,  \xi_2 > 0$, $\alpha + \beta = 2^*_s$ and that one of the following conditions hold,
\begin{eqnarray}
 0<\mu_1 \leq\mu_2< \lambda_{1,s}, \;  \label{eq.1.50}
\end{eqnarray}
\begin{eqnarray}
\lambda_{k,s}< \mu_1 \leq \mu_2< \lambda_{k+1,s}, \; for \; some \;integer \; k \geq 0.  \label{eq.1.5}
\end{eqnarray}
Then, system (\ref{c1.1}) has a second solution.
\end{theorem}

\begin{remark}
It is important to note that the hypothesis (\ref{eq.1.50}) implies that the conditions (\ref{eq:1.3}) and (\ref{eq:1.4}) are verified and the hypothesis (\ref{eq.1.5}) implies in
(\ref{eq:1.3}) and (\ref{eq:1.40}). In both cases, $\det({\lambda_{1,s}}\; I-A) > 0$.
\end{remark}

\begin{theorem}\label{teo4}
Suppose $N>6s$ and
\begin{eqnarray} \nonumber
\xi_1, \xi_2 > 0 \;and \;\lambda_{k,s}= \mu_1 \leq \mu_2<\lambda_{k+1,s}, \; for \; some \; k >1.
\end{eqnarray} In addition assume that
\begin{eqnarray} F_1 = ( f_1 , g_1 ) \in (Ker ((-\overrightarrow{\Delta})^s -\lambda_{k,s} I))^\bot.
\label{eq:1.400}
\end{eqnarray}
Then system (\ref{c1.1}) has a second solution.
\end{theorem}

\section{Notations and preliminary stuff}
For any measurable function $u:\R^N\to\R$ the Gagliardo seminorm is defined by 
\[
[u]_{s}:=\Big(C(N,s)\int_{\R^{2N}}\frac{|u(x)-u(y)|^2}{|x-y|^{N+2s}} dxdy\Big)^{1/2}=
\Big(\int_{\R^{N}}|(-\Delta)^{\frac{s}{2}} u|^2dx\Big)^{1/2}.
\]
The second equality follows by \cite[Proposition 3.6]{nezza} when 
the above integrals are finite. Then, we consider the fractional Sobolev space
\[
H^s(\R^N)=\{u\in L^2(\R^N):\,\, [u]_{s}<\infty\},\quad
\|u\|_{H^s}=(\|u\|^2_{L^2}+[u]_s^2)^{1/2},
\]
which is a Hilbert space. We use the closed subspace
\begin{equation*}
X_0^s(\Omega):=\{u\in H^s(\R^N):\,u=0\,\,\,\mbox{a.e. in $\R^N\setminus\Omega$}\}.
\end{equation*}
By Theorems 6.5 and 7.1 in \cite{nezza}, the imbedding 
$X_0^s(\Omega)\hookrightarrow L^r(\Omega)$
is continuous for $r \in [1, 2^*_s ]$ and
compact for  $r \in [1, 2^*_s )$.
Due to the fractional Sobolev inequality, $X_0^s(\Omega)$ is a Hilbert space with inner product
\begin{equation*}
\langle u,v\rangle_{X_0^s}:=C(N,s)\int_{\R^{2N}}\frac{(u(x)-u(y))(v(x)-v(y))}{|x-y|^{N+2s}}\, dx\, dy,
\end{equation*}
which induces the norm $\|\cdot\|_{X_0^s}=[\,\cdot\,]_{s}$.
Observe that by Proposition 3.6 in \cite{nezza}, we have the following identity
$$\| u \|^2_{X_0^s} = \frac{2}{C(N,s)} \| (-\Delta)^{\frac{s}{2}} u \|^2_{\R^N}, \; u \in X_0^s(\Omega).$$
Then it is proved that for $u,v \in X_0^s(\Omega),$
\begin{equation*}
\frac{2}{C(N,s)}\int_{\R^{N}} u(x) (-\Delta)^s v(x) dx=\int_{\R^{2N}}\frac{(u(x)-u(y))(v(x)-v(y))}{|x-y|^{N+2s}}\, dx\, dy,
\end{equation*}
in particular, $(-\Delta)^s$ is self-adjoint in $X_0^s(\Omega).$

Now, we consider the Hilbert space given by the product space
\begin{equation}\nonumber
Y(\Omega):=X_0^s(\Omega)\times X_0^s(\Omega),
\end{equation}
equipped with the inner product
\begin{equation}\nonumber
\langle (u,v),(\varphi,\psi)\rangle_Y:=\langle u,\varphi\rangle_{X_0^s}+\langle v,\psi\rangle_{X_0^s}
\end{equation}
and the norm
\begin{equation}\nonumber
\|(u,v)\|_{Y}:=(\|u\|_{X_0^s}^2+\|v\|_{X_0^s}^2)^{1/2}.
\end{equation}
The space $L^{r}(\Omega) \times L^{r}(\Omega)$ ($r>1$) is considered with the standard product norm
\begin{equation}\nonumber
\|(u,v)\|_{L^{r} \times L^{r}}:=(\|u\|_{L^{r}}^2+\|v\|_{L^{r}}^2)^{1/2}.
\end{equation}
Besides, we recall that 
\begin{equation}
\label{controllo}
\mu_1 |U|^2\leq (AU,U)_{\R^2}\leq \mu_2 |U|^2,\quad\text{for all $U:=(u,v) \in\R^2$,}
\end{equation}
where $\mu_1 \leq \mu_2$ are the eigenvalues of the symmetric matrix $A.$
In this paper, we consider the following notation for product space $ \S \times \S := \S^2$
and $$ \qquad w^+ (x):=\max \{w(x),0\},\,\, w^- (x):=\max\{-w(x),0\}
$$ for positive and negative part of a function $w.$  Consequently we get $w=w^+ - w^- .$

Since we are wanted to obtain a solution for the problem \eqref{eq:1.1} with critical growth, we defined $S$ be the best constant for the Sobolev-Hardy embedding
$$X_0^s(\Omega)\hookrightarrow L^{2^*_s}(\Omega).$$

The constant
\begin{align*}
S=S_{\alpha+\beta}(\Omega)=\inf_{u\in X_0^s(\Omega)\backslash \{0\}}\left\{\dfrac{ \| u \|^2_{X_0^s}}{\left( \displaystyle\int_{\Omega}\lvert u\rvert^{2^*_s}dx\right) ^{2/2^*_s}}\right\}.
\end{align*}

In \cite{fm},
Chen, Li and Ou prove that the best Sobolev constant $S_{\alpha+\beta}=S$ is achieved by $w,$ where $w$ is the unique positive solution (up to translations and dilations) of
\begin{align*}
(-\Delta)^s w=w^{2^*_s-1},\quad \mbox{in}\;\;\R^{N},\quad w\in L^{2^*_s}(\Omega).
\end{align*}

For the case of problems involving systems, we need the following definition.

\begin{align*}
S_s=S_s(\alpha,\beta)(\Omega)=\inf_{(u,v)\in  Y\backslash \{0\}}\dfrac{\| (u,v)\|_Y^2}{\left( \displaystyle\int_{\Omega}\lvert u\rvert^{\alpha}\lvert v\rvert^{\beta}+\xi_1|u|^{\alpha+\beta}+\xi_2|v|^{\alpha+\beta}dx\right) ^{2/2_s^*}}.
\end{align*}

The following result establishes a relationship between $S$ and $S_s$. In local case, it was proved in \cite{z}, which the proof in our case follows arguing as was done there combined with the arguments in \cite{f} and \cite{lfo} for the nonlocal case.

\begin{lemma}\label{relation} Let $\Omega$ be a domain
(not necessarily bounded), then there exists a positive constant $m$ such that $S_s=mS$.
Moreover, if $w_0$ achieves $S$ then $(s_0 w_0, t_0 w_0)$ achives $S_s$ for some positive constants $s_0$ and $t_0$.
\end{lemma}

\begin{remark} \label{OBS10}
The constant $m$ of the previous lemma is given by $m=M^{-1}$, where
$M=\max J(s,t)$ is attained in some $(B , C )$
(with $ B , C > 0$) of the compact set $\{(s, t)  \in \mathbb{R}^2 \; : \; |s|^2 + |t|^2 = 1 \}$
with $$J(s,t):= (|s|^\alpha |t|^\beta + \xi_1 |s|^{\alpha+\beta} + \xi_2 |t|^{\alpha+\beta})^{\frac{2}{\alpha+\beta}}.$$ Therefore,
$$\frac{B^2 + C^2}{(B^\alpha C^\beta + \xi_1 B^{\alpha+\beta} + \xi_2 C^{\alpha+\beta})^{\frac{2}{\alpha+\beta}}} =m.$$
\end{remark}

\subsection{An eigenvalue problem}

For $\lambda \in \R$, we consider the problem with homogeneous Dirichlet boundary condition

\begin{equation} \label{c1.300}
\left\{
\begin{aligned}
(-\Delta)^s u &= \lambda u  && \text{in $\Omega$,} \\
u =& \; 0 && \text{in $\R^N\setminus\Omega$}.
\end{aligned}
\right.
\end{equation}
If (\ref{c1.300}) admits a weak solution $u \in X_0^s(\Omega) \setminus \{0 \}$, then
$\lambda$ is called an eigenvalue and $u$ a $\lambda$-eigenfunction. The set of all eigenvalues is referred as
the spectrum of $(-\Delta)^s$ in $X_0^s(\Omega)$ and denoted by $\sigma((-\Delta)^s).$
Since $K=[(-\Delta)^s]^{-1}$ is a compact operator, the problem (\ref{c1.300}) can be written as $u=\lambda Ku$ with $u \in L^2(\Omega),$ hence the following results are true (see  \cite{servadei}, \cite{servadeiE}).

$(i)$ problem (\ref{c1.300}) admits an eigenvalue $\lambda_{1,s}= \min \sigma((-\Delta)^s)>0$ that can be characterized
as follows
\begin{equation}\label{eigen1}
\lambda_{1,s}=\min_{u \in X_0^s \setminus\{0\}}\frac{\displaystyle\int_{\R^{N}} |(-\Delta)^{\frac{s}{2}} u (x)|^2 dx}{\displaystyle\int_{\R^N}|u(x)|^2 dx};
\end{equation}
		
$(ii)$ there exists a non-negative function $\varphi_{1,s} \in X_0^s(\Omega)$, which is an eigenfunction corresponding to $\lambda_{1,s}$, attaining the minimum in (\ref{eigen1});

$(iii)$ all $\lambda_{1,s}$-eigenfunctions are proportional, and if $u$ is a $\lambda_{1,s}$-eigenfunction, then either $u(x) > 0$ a.e. in $\Omega$ or $ u(x) < 0$ a.e. in $\Omega$;

$(iv)$ the set of the eigenvalues of problem (\ref{c1.300}) consists of a sequence $\{\lambda_{k,s}\}$ satisfying
$$0< \lambda_{1,s} < \lambda_{2,s}\leq \lambda_{3,s}\leq \ldots \leq \lambda_{j,s}\leq \lambda_{j+1,s}\leq  \ldots, \ \lambda_{k,s}\rightarrow \infty,\ \mbox{as}\ k \rightarrow \infty,$$
which is characterized by
\begin{equation}\label{eigenvalue}
\lambda_{k+1,s}=\min_{u \in \mathbb{P}_{k+1}\setminus\{0\}}\frac{\displaystyle\int_{\R^{2N}}\frac{|u(x)-u(y)|^2}{|x-y|^{N+2s}} dxdy}{\displaystyle\int_{\R^N}|u(x)|^2 dx}
\end{equation}
where
\begin{equation*}
\mathbb{P}_{k+1}=\{ u \in X_0^s(\Omega):  \ \; \langle u, \varphi_{j,s} \rangle_X =0, \ j=1,2,\ldots, k\};
\end{equation*}
$(v)$ if $\lambda \in \sigma((-\Delta)^s) \setminus \{\lambda_{1,s} \} $ and $u$ is a $\lambda$-eigenfunction, then $u$ changes sign in $\Omega.$\\
$(vi)$ Denote by $\varphi_{k,s}$ the eigenfunction associated to the eigenvalue $\lambda_{k,s},$ for each $k\in \N.$
The sequence $\{\varphi_{k,s}\}$  is an orthonormal basis either of $L^2(\Omega)$ or of $X_0^s(\Omega)$.
		
\begin{remark}
Every eigenfunction of $(-\Delta)^s $ is in  $C^{0, \sigma} (\overline{\Omega})$ for some $\sigma \in (0, 1)$ (see Theorem 1 of \cite{servadei} or Proposition 2.4 of \cite{servadeiY}).
\end{remark}

\section{Proof of Theorem \ref{teo1}}

The proof of the Theorem \ref{teo1} needs the following lemma (see details in \cite{frp}).

\begin{lemma} \label{lema2.1}
If \eqref{eq:1.3} hold and $F_1 \in
L^2(\Omega)\times L^2(\Omega)$, then the system
\begin{equation}\label{pp2}
\left\{
\begin{aligned}
(-\overrightarrow{\Delta})^s U &= AU + F_1 &&  \text{in $\Omega$,} \\
U = & \;0  && \text{in $\R^N\setminus\Omega$},
\end{aligned}
\right.
\end{equation}
has a unique solution $U_0=(u_0,v_0) \in Y(\Omega).$
\end{lemma}

\begin{remark}
If (\ref{eq:1.400}) holds, using the Fredholm alternative, we have that (\ref{pp2}) has a unique solution.
\end{remark}
\begin{remark}\label{obs1}
If $F_1 \in L^q (\Omega)\times L^q (\Omega)$ with $q> \frac{N}{2s}$, by [\cite{le}, Theorem 3.13], we know that the solution $U_0=(u_0,v_0) \in  C^{0}(\overline{\Omega})\times C^{0}(\overline{\Omega}).$\\
If $F_1 \in L^\infty (\Omega)\times L^\infty (\Omega)$, by [\cite{xj}, Proposition 4.6], the solution $U_0=(u_0,v_0) \in  C^{0,s}(\overline{\Omega})\times C^{0,s}(\overline{\Omega}).$\\
If $s=1/2$ and $F_1 \in C^{0,\sigma}_0(\overline{\Omega}) \times C^{0,\sigma}_0(\overline{\Omega})$, with $0<\sigma <1$ and $N > 2s,$ then
$U_0 \in C^{1,\sigma}(\overline{\Omega}) \times C^{1,\sigma}(\overline{\Omega})$ and
$\| U_0 \|_{(C^{1,\sigma}(\overline{\Omega}))^2} \leq c \| F_1 \|_{(C^{0,\sigma}(\overline{\Omega}))^2}$ (see \cite{cata} Proposition 3.1) and
if $s> 1/2,$ arguing as in \cite{bcps}, we have that $U_0 \in C^{1,2s-1}(\overline{\Omega}) \times C^{1,2s-1}(\overline{\Omega}).$
Moreover, a bootstrap argument ensures that if the
function $F_1 \in C^0(\overline{\Omega}) \times C^0(\overline{\Omega})$ and $N > 2s$, then the solution
$U_0$ given by Lemma \ref{lema2.1} satisfies $\| U_0 \|_{(C^{0,\sigma}(\R^N))^2} \leq c \| F_1 \|_{(L^q(\Omega))^2},$ where $\sigma= min \{s, 2s - \frac{N}{q}\},$
for some constant depending only on $N,s, q$ and $\Omega$ (see \cite{ros} Proposition 1.4).
\end{remark}

We are ready to prove the existence of a negative solution for the system \eqref{c1.1}.

{\bf Proof of Theorem \ref{teo1}.} We will prove the theorem when the conditions \eqref{eq:1.3} and \eqref{eq:1.40} hold
(other cases (\eqref{eq:1.3} and (\ref{eq:1.4}) or (\ref{eq:1.400})) are analogous to this and left to the reader).

By Lemma \ref{lema2.1} and Remark \ref{obs1}, the system
\begin{equation*} 
\left\{
\begin{aligned}
(-\overrightarrow{\Delta})^s U &= AU + F_1 && \text{in $\Omega$,} \\
U = & \;0  && \text{in $\R^N\setminus\Omega$},
\end{aligned}
\right.
\end{equation*}
has a unique solution $ U_0=(u_0 ,  v_0) \in C^{0}(\overline{\Omega}) \times C^{0}(\overline{\Omega}).$\\
Besides, $(w,z)= \left(\dfrac{(\lambda_{1,s} - c)t + br}{\det({\lambda_{1,s}}I-A)} \phi_{1,s} ,  \dfrac{bt + (\lambda_{1,s} - a)r}{\det({\lambda_{1,s}}I-A)}
\phi_{1,s} \right)$ is the unique solution of the system
\begin{equation*} 
\left\{
\begin{aligned}
(-\overrightarrow{\Delta})^s U &= AU + T \phi_{1,s} && \text{in $\Omega$,} \\
U = & \;0  && \text{in $\R^N\setminus\Omega$}.
\end{aligned}
\right.
\end{equation*}
Consequently, if
\begin{eqnarray*} && u_{T}= \frac{(\lambda_{1,s} - c)t + br}{\det({\lambda_{1,s}}I-A)} \phi_{1,s} + u_0 ,\\
&& v_{T} = \frac{bt + (\lambda_{1,s} - a)r}{\det({\lambda_{1,s}}I-A)}
\phi_{1,s} + v_0,
\end{eqnarray*}
then $U_T =(u_T, v_T)$ is a solution of the system
\begin{equation*}
\left\{
\begin{aligned}
(-\overrightarrow{\Delta})^s U &= AU + T \phi_{1,s} + F_1 && \text{in $\Omega$,} \\
U = & \;0   && \text{in $\R^N\setminus\Omega$}.
\end{aligned}
\right.
\end{equation*}

Clearly if  $ u_{T}$ and $v_{T}$ are negative in $\Omega$, we deduce also that $U_T$ is a solution of (\ref{c1.1}).
Therefore, to conclude the proof under the conditions \eqref{eq:1.3} and \eqref{eq:1.40} (see Remark \ref{obs0}), it suffices to show the existence of an unbounded
region $\textbf{\textit{R}} \subset \mathbb{R}^2$ where $ u_{T}$ and $ v_{T}$ are negative in $\Omega$ for every $ T=(t , r) \in \textbf{\textit{R}}.$

Indeed, since $ \phi_{1,s} \in C^{0,\sigma}(\overline{\Omega})$ is strictly positive in $\Omega$ (see corollary 4.8 in \cite{moli}) and $u_0, v_0 \in C^{0}(\overline{\Omega})$, there exists
$\eta, \vartheta \ll 0$ such that
$$ \eta \varphi_{1,s} + u_0  < 0 \; \mbox{in} \;\Omega ,$$
$$ \vartheta \varphi_{1,s} + v_0  < 0 \; \mbox{in} \;\Omega.$$

Then  $ u_{T}$ and $ v_{T}$ are negative in $\Omega$ for every $ T=(t , r) \in \textbf{\textit{R}}$ and the proof of theorem is concluded. $\hfill \rule{2mm}{2mm}$

\section{Proof of Theorem \ref{teo2}}

Let $U_T:= (u_T, v_T)$ be the negative solution with $u_T, v_T <0$ in $\Omega$ given by Theorem \ref{teo1} for $T \in \textbf{\textit{R}}.$ Notice that if $ \overline{U} \neq (0,0)$ is a solution of
\begin{equation}\label{c1.7}
\left\{
\begin{aligned}
(-\overrightarrow{\Delta})^s U &= AU + \nabla F(U+U_T) && \text{in $\Omega$,} \\
U = & \;0  && \text{in $\R^N\setminus\Omega$},
\end{aligned}
\right.
\end{equation}
then $ U = \overline{U} + U_T $  is a (second) solution of the system (\ref{c1.1}) with $\overline{U} + U_T \neq U_T. $ Therefore, to prove the Theorem \ref{teo2}, we only have to show that the system (\ref{c1.7}) has a nonzero solution for every $T \in \textbf{\textit{R}}.$

Observe that the weak solutions of (\ref{c1.7}) are the critical
points of the functional $ \I_{\lambda,s}: Y(\Omega) \longrightarrow \mathbb{R}$ given by
\begin{align*}
\I_{\lambda,s}(U) =& \frac{C(N,s)}{2} \displaystyle\int_{\R^{2N}}\frac{|u(x)-u(y)|^2+|v(x)-v(y)|^2}{|x-y|^{N+2s}} dxdy\\
-& \frac{1}{2} \displaystyle\int_{\Omega} (AU,U)_{\R^2} dx -\int_{\Omega}  F(U+U_T) dx ,
\end{align*}
where $$F(U) := \frac{1}{\alpha + \beta} \left[{u}^{\alpha}_{+}
{v}^{\beta}_{+}  + \xi_1 {u_+}^{\alpha+\beta} + \xi_2 {v_+ }^{\alpha+\beta} \right], \; \mbox{for every} \; U=(u,v) \in \R^2$$
and that $U=0$ is a critical point for $\I_{\lambda,s}$ with $\I_{\lambda,s}(0)=0.$
\begin{remark}\label{obs10} The nonlinearity $F$ is $(\alpha+\beta)$-homogeneous, i.e.
$$F(\lambda U) = \lambda^{\alpha+\beta} F(U), \; \forall U \in \R^2, \; \forall \lambda \geq 0.$$
In particular:
\begin{itemize}
\item[\textit{\textbf{(i)}}]	  $ (\nabla F(U), U)_{\R^2}= u F_u (U) + v F_v(U) = (\alpha+\beta) F(U), \; \forall U=(u,v) \in \R^2.$
\item[\textit{\textbf{(ii)}}] $F_u$ and $F_v$ are $(\alpha+\beta-1)$-homogeneous.
\item[\textit{\textbf{(iii)}}] There exists $K>0$ such that
$$F_u(U)\leq K (|u|^{\alpha+\beta-1}+|v|^{\alpha+\beta-1}) \; \mbox{and}$$
$$F_v(U)\leq K (|u|^{\alpha+\beta-1}+|v|^{\alpha+\beta-1}),$$
\end{itemize}
for all $U=(u,v) \in \R^2.$

Since $F(U)=F(u_+ , v_+), \; \forall \; U=(u,v) \in \R^2,$ we deduce that
$$|\nabla F(U)|\leq K (u_+^{\alpha+\beta-1}+ v_+^{\alpha+\beta-1})$$ for some constant $K>0.$
\end{remark}

\subsection{Geometry of the functional $\I_{\lambda,s}$}
In this subsection, we demonstrate that the functional $\I_{\lambda,s}$ satisfies the geometric structure required by the Linking Theorem (see \cite[Theorem 5.3]{rabino}) when $ \lambda_{k,s} \leq \mu_1 \leq\mu_2< \lambda_{k+1,s}$, for some $k \geq 1$. In particular, if $\mu_2 < \lambda_{1,s}$ holds, then the functional satisfies the  conditions of the Mountain Pass Theorem.

Since $Y(\Omega)$ is a Hilbert space,
consider the following orthogonal decomposition
$Y(\Omega)= E_{k}^- \oplus E_k^{+},$ where
$$E_k^-=span\{ (0,\varphi_{1,s}),(\varphi_{1,s},0), (0,\varphi_{2,s}),(\varphi_{2,s},0), \ldots, (0, \varphi_{k,s}),(\varphi_{k,s},0) \}$$
and $ E^+_k=(E_k^-)^{\bot},$  for $1 \leq k \in \N.$ Note that $E^+_k=(\mathbb{P})^2$ and $U\in Y(\Omega),$ then $U=U^-+U^+$ with $U^-\in E^-_k$ and $U^+\in E^+_k.$

Therefore from the variational characterization (\ref{eigenvalue}), we have the following estimates:
$$ \| U\|^2_{Y} \geq \lambda_{k+1,s} \|U\|^2_{L^2 \times L^2} , \; \mbox{for all} \; U \in E^+_k,$$
$$ \| U\|^2_{Y} \leq \lambda_{k,s} \|U\|^2_{L^2 \times L^2} , \; \mbox{for all} \; U \in E^-_k.$$

Let $$S_{\rho}:=\partial B_{\rho}\cap E_k^{+}$$ and\\ $Q:=\{U\in Y(\Omega)\;:\; U=W+\zeta E,\quad W\in E_k^-,\quad \|W\|_{Y}\leq r,\quad 0\leq\zeta\leq R    \},$ where $E\in E_k^{+},$ $0<\rho<R$ and $r>0$ will be chosen later so that the following conditions hold:
\begin{align*}
\displaystyle\inf_{U\in S_{\rho}}\I_{\lambda,s}(U)\geq \sigma>0,\\
\displaystyle\max_{U\in \partial Q}\I_{\lambda,s}(U)\leq \alpha_0,\;\mbox{with}\;\alpha_0<\sigma,\\
\displaystyle\max_{U\in Q}\I_{\lambda,s}(U)\leq \dfrac{s}{N}S^{\frac{N}{2s}}.
\end{align*}

\begin{proposition}\label{GMPG}

Suppose $\Omega$ is a smooth bounded domain of $\R^N$, $\alpha+\beta=2^{*}_{s}$ and $\lambda_{k,s} < \mu_1 \leq\mu_2< \lambda_{k+1,s},$ for some $k \in \N.$ Then
there exists $\rho_0>0$ and a function $\alpha:[0,\rho_0]\rightarrow \R^+$ such that $$\I_{\lambda,s}(U)\geq \alpha(\rho)\;\;\mbox{for all}\;\; U \in S_{\rho}:=\partial B_\rho(0) \cap E_k^{+}.$$
\end{proposition}
Explicitly the maximum value of $\alpha(\rho)$ is
\begin{align}\label{alphachapeu}
\hat{\alpha}=\frac{s}{N}S^{N/2s}(1-\frac{\mu_2}{\lambda_{k+1,s}})^{\frac{N}{2s}}\frac{1}{(1+\xi)^{\frac{N-2s}{2s}}}
\end{align}
and this is assumed in $\hat{\rho}=S^{\frac{N}{4s}}(1-\frac{\mu_2}{\lambda_{k+1,s}})^{\frac{N-2s}{4s}}\frac{1}{(1+\xi)^{\frac{N-2s}{4s}}}$, where $S$ is the best constant for the embedding of $X_0^s$ in $L^{2^*_s}$ and $\xi=:\max\{\xi_1,\xi_2\}$.\\

\textbf{Proof}
Using the fact that $(A(U),U)_{\R^2}\leq \mu_2|u|^2$, we obtain
\begin{eqnarray*}
\I_{\lambda,s}(U)
&\geq& \frac{1}{2}||U||^2_Y-\frac{\mu_2}{2}\int_{\Omega}|U|^2 dx-\frac{
1}{\alpha+\beta} \displaystyle\int_{\Omega}\left[\xi_1 (u+u_{r,t})_+^{2^*_s}+\right. \\
&&\left.+\xi_2 (v+v_{r,t})_+^{2^*_s}+(u+u_{r,t})^{\alpha}_+(v+v_{r,t})^{\beta}_+\right]dx.
\end{eqnarray*}
Note that
\begin{equation}\label{eq1}
s^{\alpha} t^{\beta}\leq s^{\alpha+\beta}+t^{\alpha+\beta} \; \mbox{for all} \; s,t \geq 0 \end{equation}
and
\begin{equation}\label{eq2}
\int_{\Omega}(u+u_{r,t})^{2^*_s}_+ dx\leq \int_{\Omega}|u|^{2^*_s} dx \leq S^{-2^*_s/2}||u||^{2^*_s}_{X_0^s}=S^{-\frac{N}{N-2s}}||u||^{2^*_s}_{X_0^s}.
\end{equation}
Similarly
\begin{equation}\label{eq3}
\int_{\Omega}(v+v_{r,t})^{2^*_s}_+ dx \leq S^{- \frac{N}{N-2s}} ||v ||^{2^*_s}_{X_0^s}.
\end{equation}
Then, by (\ref{eq1}), (\ref{eq2}) and (\ref{eq3}), we have
\begin{align*}
\I_{\lambda,s}(U) &\geq   \frac{1}{2}\left(1-\frac{\mu_2}{\lambda_{k+1,s}}\right)||U||^2_Y\\
&-\left(\frac{
(1+\xi_1)}{\alpha+\beta} S^{-\frac{N}{N-2s}}||u||^{2^*_s}_{X_0^s}+\frac{(1+\xi_2)}{\alpha+\beta}S^{-\frac{N}{N-2s}}||v||^{2^*_s}_{X_0^s}\right).
\end{align*}
Since $\xi=:\max\{\xi_1,\xi_2\}\geq \xi_1,\xi_2$, we obtain

\begin{equation*}
\I_{\lambda,s}(U) \geq\frac{1}{2} \left(1-\frac{\mu_2}{\lambda_{k+1,s}} \right){\rho}^2-\frac{
(1+\xi)}{\alpha+\beta} S^{-\frac{N}{N-2s}}{\rho}^{2^*_s}=:\alpha(\rho),
\end{equation*}
where $\rho=||U||_Y$.
Using a standart calculus argument, we obtain that the maximum of $\alpha(\rho)$ is attained in
\begin{equation*}
\rho_0=\frac{1}{(1+\xi)^\frac{N-2s}{4s}}S^{N/4s}\left( 1-\frac{\mu_2}{\lambda_{k+1,s}}\right)^{\frac{N-2s}{4s}}.
\end{equation*}
So, the function $\alpha:[0,\rho_0]\rightarrow \R^+$ is such that $\I_{\lambda,s}(U)\geq\alpha(\rho)$ for all $U\in S_p$ and the maximum value is
\begin{equation}\label{rho.chapeu}
\alpha(	\rho_0)=\frac{s}{N}S^{N/2s}\left(1-\frac{\mu_2}{\lambda_{k+1,s}}\right)^{\frac{N}{2s}}\frac{1}{(1+\xi)^{\frac{N-2s}{2s}}}.
\end{equation}
Therefore, $\I_{\lambda,s}(U)\geq \alpha(\rho)\;\;\mbox{for all}\;\; U \in S_{\rho}.$ The proof of the proposition is complete.
$\hfill \rule{2mm}{2mm}$


It is well know (see \cite[Theorem 1.1]{fmpsz}) that $S=S_{\alpha+\beta}$ is achieved by
\begin{align}\label{u.til}
\widetilde{u}(x):=k(\mu^2+| x-x_0 |^2)^{-\frac{N-2s}{2}}, \end{align}
with $k\in\R\setminus\{0\}, \mu>0$ and $x_0\in\R^N$ fixed constants.\\

Equivalently, we see that
\begin{align*}
S=\inf\limits_{\begin{tiny} \begin{array}{cc}
u \in X_0^s\setminus\{0\} \\
\|u\|_{L^{2^{*}_s}}=1\\
\end{array} \end{tiny}} \int_{\R^{2N}}\frac{|u(x)-u(y)|^2}{|x-y|^{N+2s}} dxdy = \int_{\R^{2N}} \frac{|\overline{u}(x)- \overline{u}(y)|^2}{|x-y|^{N+2s}}dxdy
\end{align*}
where $\overline{u}(x)=
\dfrac{\widetilde{u}(x)}{||\widetilde{u}||_{L^{2^*_s}}}.
$
By translation, supose $x_0=0$ in (\ref{u.til}).
Then, the function $u^*(x)=\overline{u}\left(\frac{x}{S^{\frac{1}{2s}}}\right)$, $x\in \R^N$, is a solution for the problem
\begin{align}\label{pbm20}
(-\Delta)^su=|u|^{2^*_s-2},\;\mbox{ em } \R^N
\end{align}
satisfying
\begin{align*}
||u^*||^{2^*_s}_{L^{2^*_s}(\R^N)}=S^{\frac{N}{2s}}.
\end{align*}
As in \cite{servadeiTAMS}, for every $\epsilon>0$ we define the family of functions
$$U_{\epsilon}(x):=\epsilon^{-\frac{N-2s}{2}}u^{*}\left(\dfrac{x}{\epsilon}\right),\; x\in \R^N,$$
then $U_{\epsilon}$ is a solution of (\ref{pbm20}) and verify for all $\epsilon>0$
\begin{align*}
\int_{\R^{2N}}\frac{|U_{\epsilon}(x)-U_{\epsilon}(y)|^2}{|x-y|^{N+2s}} dxdy = \int_{\R^{2N}} |U_{\epsilon}(x)|^{2^*_s}dxdy = S^{\frac{N}{2s}}.
\end{align*}
Now, take a fixed $\delta>0$ such that $B_{4\delta}\subset\Omega.$ Let $\eta\in C_{c}^{\infty}(\R^{N})$ be a cut-off function such that $0\leq\eta\leq 1$ in $\R^N,$ $\eta=1$ in $B_{\delta}$ and $\eta=0$ in $\R^N\setminus B_{2\delta}$, where $B_r=B_r(0)$ is the ball centered in origin and with radius $r>0$.\\
Define the family of nonnegative truncated functions
\begin{align}\label{u.epsilon}
u_{\epsilon}(x):=\eta(x)U_{\epsilon}(x)\;\; x\in \R^N,
\end{align}
and note that $u_{\epsilon}\in X^s_0$.

The following Brezis-Nirenberg estimates for nonlocal setting was proved in \cite{servadeiTAMS} (also see \cite{low}), which are similar to those proved for the local case in \cite{bn}.
\begin{lemma}\label{prop2} Suppose $s \in (0,1)$ and $N>2s,$ then for $\epsilon>0$ small enough, the following estimates hold true,
\begin{align}
\int_{\R^{2N}}\frac{|u_\epsilon(x)-u_\epsilon(y)|^2}{|x-y|^{N+2s}} dxdy \leq S^{N/2s}+ O(\epsilon^{N-2s}),\nonumber
\end{align}
\begin{align}
\nonumber \int_{\R^N} |u_\epsilon(x)|^{2} dx \geq
\left\{\begin{array}{lcr} C_s \epsilon^{2s} + O(\epsilon^{N-2s}) &\mbox{if}&N > 4s,\\
\nonumber C_s \epsilon^{2s}|log \epsilon| + O(\epsilon^{2s}) &\mbox{if}&N = 4s,\\
\nonumber C_s \epsilon^{N-2s} + O(\epsilon^{2s}) &\mbox{if}&2s< N\leq  4s,\\
\end{array}\right.
\end{align}
\begin{align}
\int_{\R^N} |u_\epsilon(x)|^{2^{*}_s} dx=S^{N/2s}+ O(\epsilon^{N}),\nonumber
\end{align}
\begin{align}
||u_\epsilon||_{L^{1}(\R^N)}= O(\epsilon^{\frac{N-2s}{2} }),\nonumber
\end{align}
\begin{align}
||u_\epsilon||^{2^*_s-1}_{L^{2^*_s-1}(\R^N)}= O(\epsilon^{\frac{N-2s}{2} }).\nonumber
\end{align}
\end{lemma}

Denote by $P_-$ the ortogonal projection of $X_0^s$ in $B^-_k=span\{\phi_1,\phi_2,...,\phi_k\}$ and $P_+$ the orthogonal projection of $X_0^s$ in $A^+_k:=(B^-_k)^\bot$.\\
Chosing depending on $\epsilon>0$ the vetorial function given by
\begin{align*}
e=\vec{e}_{\epsilon}=(B(P_+u_{\epsilon}),C(P_+u_{\epsilon})) \in E^+_k,
\end{align*} where $u_\epsilon$ is given in $\eqref{u.epsilon}$ and $B$ and $C$ are given by Remark \ref{OBS10}.

We will denote $P_+u_{\epsilon}$ by $e_{\epsilon}$ and consequently $\vec{e}_{\epsilon}=(B e_{\epsilon},Ce_{\epsilon})$.
\begin{remark}
(i) $e_{\epsilon}\in A^+_k$; \\ (ii) $ \langle(B e_{\epsilon},Ce_{\epsilon}), (0, \phi_j)\rangle_{L^2\times L^2} =0=\langle(B e_{\epsilon},C e_{\epsilon}), (\phi_j, 0)\rangle_{L^2\times L^2} $, for all $  j=1, ..., k.$ Then $e=\vec{e}_{\epsilon} \in E^+_k$.
\end{remark}
Hence, the following results was proved in \cite{AMBROSIO},
which are similar to those proved for the local case in \cite{djj}.

\begin{lemma}\label{lemm2} For $s \in (0,1)$ and $N>2s,$
then for $\epsilon>0$ small enough, the following estimates hold true,
	\begin{align}
||P_+u_\epsilon||^{2}_{X_0^s}\leq [u_\epsilon]_s^2 \leq S^{N/2s}+ O(\epsilon^{N-2s}),\nonumber 
\end{align}
\begin{align}
\left| ||P_+u_\epsilon||^{2^*_s}_{L^{2^*_s}(\Omega)}-||u_\epsilon||^{2^*_s}_{L^{2^*_s}(\Omega)}\right|\leq C\epsilon^{N-2s},\nonumber
\end{align}
\begin{align}
||P_+u_\epsilon||_{L^{1}(\Omega)}\leq C\epsilon^{\frac{N-2s}{2} },\nonumber
\end{align}
\begin{align}
||P_+u_\epsilon||^{2^*_s-1}_{L^{2^*_s-1}(\R^N)}\leq C\epsilon^{\frac{N-2s}{2} },\nonumber
\end{align}
\begin{align}\label{limitado}
|P_-u_\epsilon(x)|\leq C\epsilon^{\frac{N-2s}{2} },\; for \; x\in\Omega.
	\end{align}
\end{lemma}

Fix $K>0$ and define $\Omega_{\epsilon,K}=\{x\in \Omega : e_{\epsilon} (x)=(P_+u_{\epsilon})(x)>K\}.$ By \eqref{limitado} we can deduce
\begin{eqnarray*}
e_{\epsilon}(0)=(P_+u_{\epsilon})(0)=u_{\epsilon}(0)-P_-u_{\epsilon}(0)
\geq\frac{C_0}{\|\tilde{u}\|_{L^{2^*_s}(\R^N)}}\epsilon^{-\frac{(N-2s)}{2}}-C\epsilon^{\frac{N-2s}{2}},
\end{eqnarray*}
which implies that $P_+u_{\epsilon}(0)\rightarrow\infty$ as  $\epsilon \rightarrow 0$
By the continuity of $P_+u_{\epsilon}$, there exists $\nu>0$ such that $B_{\nu}\subset\Omega_{\epsilon,K}.$ Therefore, we have the result below.
\begin{lemma}\label{lemm} For $s \in (0,1)$ and $N>2s,$ we have
\begin{align}
||P_+u_\epsilon||^{2^*_s}_{L^{2^*_s}(\Omega_{\epsilon,K})}=||u_\epsilon||^{2^*_s}_{L^{2^*_s}(\Omega)}+ O(\epsilon^{N-2s}).\nonumber
\end{align}
\begin{align}
||P_+u_\epsilon||^{2^*_s-1}_{L^{2^*_s-1}(\Omega_{\epsilon,K})}=||u_\epsilon||^{2^*_s-1}_{L^{2^*_s-1}(\Omega)}+ O(\epsilon^{\frac{N+2s}{2}}).\nonumber
\end{align}
\begin{align}
||P_+u_\epsilon||_{L^{1}(\Omega_{\epsilon,K})}=||u_\epsilon||_{L^{1}(\Omega)}+  O(\epsilon^{N}).\nonumber
	\end{align}
\end{lemma}

To prove the geometric conditions of the Linking Theorem, we need two results that can be found in \cite{djj} and \cite{fd} for the case when $s=1$. The proof is similar for $s\in (0,1)$.

\begin{lemma}\label{L2} Given $u,v\in L^p(\Omega)$ with
$2\leq p\leq 2^*_s$ and $u+v>0$ a.e. on a measurable subset
$\Sigma\subset\Omega$, it holds
\begin{equation*}
\Big| \int_{\Sigma} (u+v)^p dx - \int_{\Sigma} |u|^p dx -
\int_{\Sigma} |v|^p dx \Big|\leq C \int_{\Sigma}(|u|^{p-1}|v|+
|u||v|^{p-1})\, dx,
\end{equation*}
with a constant $C>0$ depending only on $p$.
\end{lemma}

\begin{lemma}\label{L2.1} Given $(a,b), (u,v)\in L^p(\Omega) \times L^q(\Omega)$ with $p,q \geq 2$ and
$p+q\leq 2^*_s.$ If $a+b, \; u+v>0$ a.e. on a measurable subset
$\Sigma\subset\Omega$ and $H(x,y)=|x|^p |y|^q,$ then
\begin{eqnarray}\nonumber
&\Big| \displaystyle\int_{\Sigma} H(a+u, b+v) dx - \int_{\Sigma} H(u,v) dx -
\int_{\Sigma} H(a,b) dx \Big|\\ \nonumber
&\leq C \Big[ \displaystyle\int_{\Sigma}(|a|^{p-1}|b|^q |u|+
|a|^{p-1}|v|^{q}|u| + |u|^{p-1}|b|^q |a|+ |u|^{p-1}|v|^q |a|)\, dx\\ \nonumber
&+\displaystyle\int_{\Sigma}(|a|^{p-1}|v|^{q-1}|b| |u|+ |u|^{p}|b|^{q} + |u|^{p}|v|^{q-1} |b|)\, dx\\ \nonumber
&+\displaystyle\int_{\Sigma}(|a|^{p}|b|^{q-1} |v|+
|a|^{p}|v|^{q} + |u|^{p-1}|b|^{q-1} |a||v|)\, dx\\ &+\label{LP}\displaystyle\int_{\Sigma}(|b|^{q-1} |u|^p |v|+
|v|^{q-1}|a|^{p} |b|)\, dx \Big],
\end{eqnarray}
where the constant $C>0$ depending only on $p+q$.
\end{lemma}
\textbf{Proof} Let us define$$h(\zeta):=\displaystyle\int_{\Sigma} [H(a+\zeta u, b+\zeta v) - H(\zeta u, \zeta v)] \; dx.$$

Employing the Fundamental Theorem of the Calculus,
$
|h(1) - h(0)| = \displaystyle\int_0^1 h'(\zeta) \; d \zeta ,$ and consequently
\begin{eqnarray} \nonumber
&\displaystyle\int_{\Sigma} [H(a+ u, b+ v) - H(a, b)] \; dx\\ \label{TFC}
&\leq \displaystyle\int_0^1\int_{\Sigma} |((\nabla H(a+ \zeta u, b+ \zeta v) - \nabla H(\zeta u, \zeta v))\; , \; (u,v))_{\R^{2}}| \; dx d\zeta.
\end{eqnarray}

Using the Mean Value Theorem to the function $ \nabla H(x, y)$, there exist $\theta_1, \theta_2 \in (0, 1)$ such that
\begin{eqnarray}\nonumber
&\nabla H(a+ \zeta u, b+ \zeta v) - \nabla H(\zeta u, \zeta v)\\ \nonumber
&=\Big(p |a+\zeta u|^{p-2}(a+\zeta u)|b+\zeta v|^q - p |\zeta u|^{p-2}(\zeta u)|\zeta v|^q \; \; , \\ \nonumber
&q |a+\zeta u|^{p} |b+\zeta v|^{q-2} (b+\zeta v) - q |\zeta u|^{p} |\zeta v|^{q-2} (\zeta v)\Big)\\ \nonumber
&=\Big(p(p-1) |(1-\theta_1)a+\zeta u|^{p-2}|(1-\theta_1)b+\zeta v|^{q} a\\ \nonumber
&+ pq |(1-\theta_1)a+\zeta u)|^{p-2} ((1-\theta_1)a+\zeta u)|(1-\theta_1)b+\zeta v|^{q-2} (1-\theta_1)b+\zeta v)b \; , \\ \nonumber
& pq |(1-\theta_2)a+\zeta u)|^{p-2} ((1-\theta_2)a+\zeta u)|(1-\theta_2)b+\zeta v|^{q-2} (1-\theta_2)b+\zeta v)a\\ \label{TVM}
&+ q(q-1) |(1-\theta_2)b+\zeta v|^{q-2}|(1-\theta_2)a+\zeta u|^{p} b \Big).
\end{eqnarray}
Inequality \eqref{LP} follows by substituting \eqref{TVM} in \eqref{TFC} and making some forward
estimations. $\hfill \rule{2mm}{2mm}$

The following inequality which is a direct consequence of Young Inequality, is
essential for the proof of Theorem \ref{teo2}. 

\begin{lemma}\label{L2.2}
If $\alpha, \beta > 1,$
$ \alpha +\beta = 2^*_s$ and $\alpha > \dfrac{2^*_s -1}{2}$, there is $p > 2$ such that, each $\epsilon > 0,$ the following inequality holds
\begin{equation*}
|s|^\alpha |t|^\beta \leq C_{\epsilon} |s|^{2^*_s -1} + C \epsilon^p |t|^{\beta p}
\end{equation*}
where $C_\epsilon$ and $C$ are positive constants.
\end{lemma}

Finally we need

\begin{lemma}\label{L3}
Suppose $A, B, C$ and $ \theta$ positive numbers. Consider the function $\Phi_{\epsilon}(s)=\frac{1}{2}s^2A-\frac{1}{2^*_s} s^{2^*_s} B+s^{2_s^*}\epsilon^{\theta }C$ with $s>0$. Then $s_{\epsilon}=\left( \frac{A}{B-2^*_s\epsilon^{\theta}C}\right)^{\frac{1}{2^*_s-2}}$ is the maximum point of $\Phi_{\epsilon}$ and
$$\Phi_{\epsilon}(s)\leq  \Phi_{\epsilon} (s_{\epsilon}) =  \frac{s}{N} \left( \frac{A^N}{B^{N-2s}}\right)^{\frac{1}{2s}}+O(\epsilon^{\theta}).$$
\end{lemma}
\begin{lemma}\label{functional}
If $\lambda_{k,s}<\mu_1\leq \mu_2<\lambda_{k+1,s}$, there are constants $r_0,R_0>0$ and $\epsilon_0>0$ such that, for $r>r_0$, $R>R_0$ and $0<\epsilon \leq \epsilon_0$, we have
$$\I_{\lambda,s}|_{\partial Q}<\hat{\alpha},$$
with $\hat{\alpha}>0$ as in Proposition $\ref{GMPG}$.
\end{lemma}
\textbf{Proof}
Let $\partial Q=\Gamma_1 \cup \Gamma_2 \cup \Gamma_3,$ where
$$ \Gamma_1 =\overline{B_R}\cap E^-_k,$$
$$\Gamma_2 =\{ U\in Y: U=W+s\vec{e}_{\epsilon} \mbox{ with } W\in E^-_k,\; ||W||_Y=r, \; 0\leq s \leq R \},$$
$$\Gamma_3 =\{ U\in Y: U=W+R\vec{e}_{\epsilon} \mbox{ with } W\in E^-_k\cap B_r(0) \}. $$
We will show that for each $ \Gamma_i$ we have $\I_{\lambda,s}|_{\Gamma_i}<\hat{\alpha},$ for all $i=1,2,3.$\\
$(i)$ For all $U\in \Gamma_1 (\subset E^-_K)$, using (\ref{controllo}), we infer that
\begin{eqnarray*}
\I_{\lambda,s}(U)\leq
\frac{1}{2}||U||_Y^2-\frac{\mu_1}{2}\frac{1}{\lambda_{k,s}}||U||_Y^2= \frac{1}{2}\left(1-\frac{\mu_1}{\lambda_{k,s}}\right) ||U||_Y^2\leq 0.
\end{eqnarray*}
$(ii)$ Let $U\in \Gamma_2$, then $U=W+s\vec{e}_{\varepsilon}$ with $W=(w_1,w_2) \in E_k^-$ and $\vec{e}:=(B(P_+ u_{\varepsilon}), C(P_+ u_{\varepsilon}))=(B e_{\varepsilon}, C e_{\varepsilon}),$
where the positive constants $B$ and $C$ are chosen as in Remark \ref{OBS10}.

Hence
\begin{eqnarray}
\nonumber \I_{\lambda,s}(U)&\leq&
\dfrac{1}{2}\left(1-\dfrac{\mu_1}{\lambda_{k,s}}\right)||W||_Y^2+\frac{s^2}{2}(B^2 +C^2)||{e}_{\varepsilon}||_{X_0^s}^2\\ \nonumber
&-&\dfrac{1}{\alpha+\beta}\displaystyle\int_{\Omega}(w_1+sBe_{\varepsilon}+u_{rt})_+^{\alpha}(w_2+sCe_{\varepsilon}+v_{rt})_+^{\beta} dx\\ \nonumber
&-&\dfrac{\xi_1}{\alpha+\beta}\displaystyle\int_{\Omega}(w_1+sBe_{\varepsilon}+u_{rt})_+^{\alpha+\beta} dx-\dfrac{\xi_2}{\alpha+\beta}\displaystyle\int_{\Omega}(w_2+sCe_{\varepsilon}+v_{rt})_+^{\alpha+\beta} dx.
\end{eqnarray}
Consider the maximum value $\hat{\alpha}$ of the function $\alpha(\rho)$ like in (\ref{rho.chapeu}), and define
\begin{equation}\label{alpha}
s_0 := \frac{\sqrt{2\frac{s}{N}S^{N/2s}(1-\frac{\mu_2}{\lambda_{k+1,s}})^{\frac{N}{2s}}\frac{1}{(1+\xi)^{\frac{N-2}{2s}}}}}{\sqrt{\displaystyle\sup_{0<\varepsilon\leq 1}||\vec{e}_{\varepsilon}||_Y^2}}\\
=\frac{\sqrt{2\hat{\alpha}}}{\sqrt{\displaystyle\sup_{0<\varepsilon\leq 1}||\vec{e}_{\varepsilon}||_Y^2}}.
\end{equation}
In order to comply the condition
$ \I_{\lambda,s}|_{\Gamma_2}<\hat{\alpha},$ we distinguish in our analysis two cases.\\
\textbf{First case:} If $0\leq s \leq s_0$.\\
The expression of $\I_{\lambda,s}$ provides the estimate
\begin{equation}\nonumber
\I_{\lambda,s}(U)\leq\dfrac{s^2}{2} ||\vec{e}_{\varepsilon}||^2_Y \leq  \dfrac{s_0^2}{2}\displaystyle\sup_{0<\varepsilon\leq 1} ||\vec{e}_{\varepsilon}||^2_Y =\hat{\alpha}.
\end{equation}
What concludes this case.\\
\textbf{Second case:} $s>s_0$.\\
Define
\begin{equation*}
K:=  \sup \Big\{\displaystyle\left\|\frac{W+(u_{r,t},v_{r,t})}{s}\right\|_{L^{\infty}\times L^{\infty}}: s_0\leq s\leq R,
\; \|W \|_Y=r, \; W\in E^-_k \Big\},
\end{equation*}
with $K>0$ independent of $R$.\\
Then, by (\ref{u.epsilon}) and \eqref{limitado}  we have
\begin{eqnarray*}
e_{\epsilon}(0)&=&(P_+u_{\epsilon})(0)=u_{\epsilon}(0)-P_-u_{\epsilon}(0)\\
&\geq & \frac{C_0}{\|\tilde{u}\|_{L^{2^*_s}(\R^N)}}\epsilon^{-\frac{(N-2s)}{2}}-c\epsilon^{\frac{N-2s}{2}}\rightarrow +\infty,
\end{eqnarray*}
as $\epsilon \rightarrow 0$ because $N>2s$.\\
By the continuity of $e_{\epsilon}$, we have
\begin{eqnarray*}
\Omega_{\epsilon}&=&\{x\in \Omega : e_{\epsilon} (x)=(P_+u_{\epsilon})(x)>K\}\neq \emptyset
\end{eqnarray*}
for $\epsilon>0$ small enough.
Therefore, by Lemmas \ref{L2} and \ref{L2.1}, for $j=1$ and $z_{r,t}= u_{r,t}$ or $j=2$ and $z_{r,t}=v_{r,t}$, we have
\small
\begin{eqnarray}\label{des}
&&\displaystyle\int_{\Omega_{\epsilon}}\left( Be_{\epsilon}+\frac{w_j+z_{r,t}}{s}\right)^{\alpha+\beta}dx
\geq  \displaystyle\int_{\Omega_{\epsilon}}|  Be_{\epsilon}|^{2^*_s} dx+ \int_{\Omega_{\epsilon}}\left|\frac{w_j+z_{r,t}}{s}\right|^{\alpha+\beta}dx \nonumber\\
&&\!\!\!\!\!\!-C\!\! \int_{\Omega_{\epsilon}}\!\!\! \left( |  Be_{\epsilon}|^{2^*_s\!\!-\!1}\left|\frac{w_j+z_{r,t}}{s}\right|\! +\! |  Be_{\epsilon}|\left|\frac{w_j+z_{r,t}}{s}\right|^{2^*_s\!\!-\!1} \! \right) dx
\end{eqnarray}
and
\begin{eqnarray}\label{des10}
&&\displaystyle\int_{\Omega_\epsilon}\!\!\! \left( \!B e_{\epsilon} \!+\!\frac{w_1 \!+\! u_{r,t}}{s} \!\!\right)_+^{\alpha}\left( \!Ce_{\epsilon} \!+\!\frac{w_2 \!+\! v_{r,t}}{s} \!\!\right)_+^{\beta}\!\!\!dx \\ \nonumber
&&\geq  \displaystyle\int_{\Omega_{\epsilon}} B^\alpha C^\beta |e_{\epsilon}|^{\alpha+\beta} dx+ \int_{\Omega_{\epsilon}}\left|\frac{w_1+u_{r,t}}{s}\right|^{\alpha}\left|\frac{w_2+v_{r,t}}{s}\right|^{\beta}dx \\ \nonumber
&&\!\!\!\!\!\!-K\!\! \int_{\Omega_{\epsilon}}\!\!\! \Big( \left|\frac{w_1+u_{r,t}}{s}\right|^{\alpha\!\!-\!1}\left|\frac{w_2+v_{r,t}}{s}\right|^{\beta} |B e_{\epsilon}|\! +\! \left|\frac{w_1+u_{r,t}}{s}\right|^{\alpha\!\!-\!1}|C e_{\epsilon}|^{\beta}|B e_{\epsilon}|\\ \nonumber
&&+ |B e_{\epsilon}|^{\alpha-1}  \left|\frac{w_2+v_{r,t}}{s}\right|^{\beta} \left|\frac{w_1+u_{r,t}}{s}\right| + |B e_{\epsilon}|^{\alpha-1}|C e_{\epsilon}|^{\beta}\left|\frac{w_1+u_{r,t}}{s}\right| \Big)\; dx
\\ \nonumber
&&\!\!\!\!\!\!-K\!\! \int_{\Omega_{\epsilon}}\!\!\! \Big( \left|\frac{w_1+u_{r,t}}{s}\right|^{\alpha\!\!-\!1}\left|\frac{w_2+v_{r,t}}{s}\right||C e_{\epsilon}|^{\beta-1}|B e_{\epsilon}|\! +\! |B e_{\epsilon}|^{\alpha} \left|\frac{w_2+v_{r,t}}{s}\right|^{\beta}\\ \nonumber
&&+|B e_{\epsilon}|^{\alpha}|C e_{\epsilon}|^{\beta-1}  \left|\frac{w_2+v_{r,t}}{s}\right| \Big)\; dx\\
\nonumber
&&\!\!\!\!\!\!-K\!\! \int_{\Omega_{\epsilon}}\!\!\! \Big( \left|\frac{w_1+u_{r,t}}{s}\right|^{\alpha}\left|\frac{w_2+v_{r,t}}{s}\right|^{\beta-1}|C e_{\epsilon}|+ \!\left|\frac{w_1+u_{r,t}}{s}\right|^{\alpha} |C e_{\epsilon}|^{\beta}\\ \nonumber
&&+ |B e_{\epsilon}|^{\alpha-1} \left|\frac{w_2+v_{r,t}}{s}\right|^{\beta-1}   \left|\frac{w_1+u_{r,t}}{s}\right| |C e_{\epsilon}|\Big)\; dx
\\
\nonumber
&&\!\!\!\!\!\!-K\!\! \int_{\Omega_{\epsilon}}\!\!\! \Big( \left|\frac{w_2+v_{r,t}}{s}\right|^{\beta-1}|B e_{\epsilon}|^{\alpha}|C e_{\epsilon}|+ |C e_{\epsilon}|^{\beta-1}\!\left|\frac{w_1+u_{r,t}}{s}\right|^{\alpha}  \left|\frac{w_2+v_{r,t}}{s}\right| \Big)\; dx.
\end{eqnarray}

Then, using the estimates \eqref{des} and \eqref{des10}, we can see that, for $\epsilon>0$ small enough,
\begin{eqnarray*}
\nonumber \I_{\lambda,s}(U) &\leq&
\frac{1}{2}\Big(1-\frac{\mu_1}{\lambda_{k,s}}\Big) ||W||_Y^2 +
\frac{s^2}{2}(B^2+C^2)||e_{\epsilon}||_{X_0^s}^2\\ \nonumber
&&- \frac{s^{2_s^*}}{2_s^*}(B^\alpha C^\beta +\xi_1 B^{2_s^*} + \xi_2 C^{2_s^*} ) 
   ||e_{\epsilon}||^{2_s^*}_{L^{2_s^*}(\Omega_\epsilon)}\\ \nonumber
&&+K\frac{s^{2_s^*}}{2_s^*} \Big(\|e_{\epsilon}\|^{2_s^*-1}_{L^{2_s^*-1}(\Omega_\epsilon)}+\|e_{\epsilon}\|_{L^{1}(\Omega_\epsilon)}+\|e_{\epsilon}\|^{\alpha+1}_{L^{\alpha+1}(\Omega_\epsilon)} + \|e_{\epsilon}\|^{\beta+1}_{L^{\beta+1}(\Omega_\epsilon)}\\
&&+ \|e_{\epsilon}\|^{\alpha-1}_{L^{\alpha-1}(\Omega_\epsilon)}+ \|e_{\epsilon}\|^{\beta-1}_{L^{\beta-1}(\Omega_\epsilon)}+\|e_{\epsilon}\|^{\beta}_{L^{\beta}(\Omega_\epsilon)}+ \|e_{\epsilon}\|^{\alpha}_{L^{\alpha}(\Omega_\epsilon)}\Big).\nonumber\\
\end{eqnarray*}
Now, for each $j \in \{ \alpha, \beta, \alpha-1, \beta-1\}$, there exists $C_j >0$ such that $$\|e_{\epsilon}\|^{j}_{L^{j}(\Omega_\epsilon)} \leq C_j \|e_{\epsilon}\|^{2_s^*-1}_{L^{2_s^*-1}(\Omega_\epsilon)}$$ and,
by lemma \ref{L2.2}, for each
$j \in \{ \alpha+1, \beta+1\}$, there exists $K_{j} >0$ such that $$\|e_{\epsilon}\|^{j}_{L^{j}(\Omega_\epsilon)} \leq K_{j} (\|e_{\epsilon}\|^{2_s^*-1}_{L^{2_s^*-1}(\Omega_\epsilon)} +\epsilon^p),$$ with $p>2.$
Therefore, using the above estimate and the
Lemmas \ref{prop2}, \ref{lemm2} and \ref{lemm}, we get
\begin{eqnarray}
\nonumber \I_{\lambda,s}(U) \leq \frac{1}{2}\Big(1-\frac{\mu_1}{\lambda_{k,s}}\Big) \|W\|_Y^2 + \Phi_{\epsilon}(s) ,
\end{eqnarray}
where $$\Phi_{\epsilon}(s):=\dfrac{s^2}{2}(B^2 +C^2)S^{\frac{N}{2s}}-\dfrac{s^{2^*_s}}{2^*_s}(B^{\alpha}C^{\beta}+ \xi_ 1 B^{2^*_s}+ \xi_2 C^{2^*_s}) S^{\frac{N}{2s}}+Ks^{2^*_s}O(\epsilon^q)$$ with $q= \min \{ \dfrac{N-2s}{2}, p \}.$
Then, applying Lemma \ref{L3}, we obtain
\begin{eqnarray*}
\I_{\lambda,s}(U) &\leq& \frac{1}{2}\Big(1-\frac{\mu_1}{\lambda_{k,s}}\Big)r^2 + \frac{s}{N} \left( \frac{\Big[(B^2+C^2)S^{\frac{N}{2s}}\Big]^N}{\Big[(B^{\alpha}C^{\beta}+ \xi_ 1 B^{2^*_s}+ \xi_2 C^{2^*_s})S^{\frac{N}{2s}}\Big]^{N-2s}}  \right)^{\frac{1}{2s}}
+O(\epsilon^{q})\\
&=& \frac{1}{2}\Big(1-\frac{\mu_1}{\lambda_{k,s}}\Big)r^2 + \frac{s}{N} \left( \frac{(B^2+C^2)^{\frac{N}{2s}}}{(B^{\alpha}C^{\beta}+ \xi_ 1 B^{2^*_s}+ \xi_2 C^{2^*_s})^{\frac{N-2s}{2s}}} \right) S^{\frac{N}{2s}}
+O(\epsilon^{q}).
\end{eqnarray*}
Since $\lambda_{k,s}<\mu_1$ and $\varepsilon > 0$ can be made arbitrarily small, we can choose $r>0$ to be arbitrarily large in the inequality above such that $\I_{\lambda,s}(U)<0$. This leads to the conclusion stated in the proposition for $U\in\Gamma_2$.\\
$(iii)$ Let $U\in\Gamma_3$. It can be expressed by $\Gamma_3$ definition as $U=W+R\vec{e}_{\varepsilon}$ with $W\in E_k^-{\cap} B_r(0)$.
Analogously to the case $(ii)$, we get
\begin{eqnarray*}
\nonumber &&\I_{\lambda,s}(U)\leq \frac{1}{2}\Big(1-\frac{\mu_1}{
\lambda_{k,s}}\Big)\|W\|_{Y}^2 + (B^2+C^2) \frac{R^2}{2} \|
e_{\epsilon}\|^2_{X_0^s}\\
&&-\frac{B^{\alpha}C^{\beta} }{2^*_s} R^{2^*_s}
\int_{\Omega}
\Big(e_{\epsilon}+\frac{w_1+u_{r,t}}{BR}\Big)^{\alpha}_+
\Big(e_{\epsilon}+\frac{w_2+v_{r,t}}{CR}\Big)^{\beta}_+ \, dx\\
&&-\frac{\xi_1 B^{2^*_s} }{2^*_s} R^{2^*_s}
\int_{\Omega}\Big(e_{\epsilon}+\frac{w_1+u_{r,t}}{BR}\Big)^{2^*_s}_+ \, dx-\frac{\xi_2 C^{2^*_s} }{2^*_s} R^{2^*_s}
\int_{\Omega}\Big(e_{\epsilon}+\frac{w_2+v_{r,t}}{CR}\Big)^{2^*_s}_+ \, dx.
\end{eqnarray*}
Due to the boundedness of the functions $W\in E^-_k\cap B_r(0) $, $u_{r,t}$ and $v_{r,t}$, there exists $k>0$ such that $\|w_1+u_{r,t}\|_{L^{\infty}}\leq k$ and $\|w_2+v_{r,t}\|_{L^{\infty}}\leq k$.
Again, since $e_\epsilon (0)=P_+u_{\epsilon}(0)\rightarrow \infty$ as $ \epsilon \rightarrow 0$, there exists $\epsilon_0 >0$ such that for all $0<\epsilon<\epsilon_0,$ we have $e_{\epsilon}(0)> 2k$. Then, by the continuity of $e_{\epsilon}$ we can find $R_1=R_1(\epsilon)>0$ and
$\eta=\eta(\epsilon)>0$ such that $|\chi|\geq \eta$ for all  $R>R_1$, where
$$\chi:=\left\{x\in \Omega \; : \; e_\epsilon (x)+\frac{w_{1}(x)+u_{r,t}(x)}{BR}>1 \; \mbox{and} \;
e_\epsilon (x)+\frac{w_{2}(x)+v_{r,t}(x)}{CR}>1\right\}.$$

Then, we find $\epsilon_0, R_0>0$ such that for $0<\epsilon<\epsilon_0$ and $R>R_0$, we have
\begin{align}
\nonumber \I_{\lambda,s}(U)\leq
0,\;\; \mbox{ for all  } U\in \Gamma_3.
\end{align}

Indeed, let $R_0>\max\{R_1,R_2\}$, where $R_0$ is such that $\alpha R_0^2-R_0^{2^*_s}<0,$ with
$$\alpha = \frac{(B^2 +C^2)2^*_s}{2(B^\alpha C^\beta + \xi_1 B^{2^*_s}+ \xi_2 C^{2^*_s})}(\eta^{-1}\|e_{\epsilon}\|_{X_0^s}^2).$$
Then, for $\epsilon >0 $ above and $R>R_0$ we find
\begin{eqnarray*}
\nonumber \I_{\lambda,s}(U)&\leq& \frac{1}{2}\Big(1-\frac{\mu_1}{
\lambda_{k,s}}\Big)\|W\|_Y^2 + (B^2+C^2)\frac{R^2}{2} \|
e_{\epsilon}\|^2_{X_0^s(\Omega)}\\
&&-\frac{ R^{2^*_s}}{2^*_s} B^\alpha C^\beta |\chi|
-\xi_1 \frac{ R^{2^*_s}}{2^*_s} B^{2^*_s}|\chi| - \xi_2 \frac{ R^{2^*_s}}{2^*_s} C^{2^*_s}|\chi|\\
&\leq& (B^2 +C^2)\frac{R^2}{2} \|
e_{\epsilon} \|^2_{X_0^s(\Omega)} -(B^\alpha C^\beta + \xi_1 B^{2^*_s}+ \xi_2 C^{2^*_s})\frac{ R^{2^*_s}}{2^*_s} \eta < 0.
\end{eqnarray*}

This completes the proof.
$\hfill \rule{2mm}{2mm}$

\begin{lemma}\label{mini}
Let $s \in  (0,1)$, $\lambda_{k,s} \leq \mu_1  \leq \mu_2<\lambda_{k+1,s}$ and $N>6s$.
Then we have the following estimate

$$ \displaystyle\max_{ \overline{Q} }\I_{\lambda,s}<\dfrac{s}{N}S^\frac{N}{2s}$$
\end{lemma}
\textbf{Proof:} Let $\epsilon < \epsilon_0$ fixed that the linking theorem geometry holds. For $W+s \vec{e}_\epsilon \in Q$, we have
\begin{eqnarray*}
\I_{\lambda,s}(W+s\vec{e}_\epsilon)& \leq &\frac{1}{2} \Big(1- \frac{\mu_1}{\lambda_{k,s}}\Big)\|W\|_Y^2 + \frac{s^2}{2} ||\vec{e}_\epsilon||_Y^2 -\frac{\mu_1}{2} s^2||\vec{e}_\epsilon||^2_{L^2 \times L^2}\\&&- \int_{\Omega}F (w+s\vec{e}+U_T) dx.
\\
\end{eqnarray*}
Let $s_0$ be defined as in \eqref{alpha}.

First case: If $0<s\leq s_0.$ Arguing as in the proof of Lemma \ref{functional} and bearing in mind \eqref{alphachapeu}, we can see that
\begin{eqnarray}\label{SS}
\I_{\lambda,s}(W+s\vec{e}_\epsilon)\leq \frac{s^2}{2} ||\vec{e}_\epsilon||_Y ^2 \leq \frac{s_0^2}{2}\displaystyle\sup_{0<\varepsilon\leq 1}||\vec{e}_{\varepsilon}||_Y^2=\hat{\alpha}<\dfrac{s}{N}\dfrac{1}{(1+\xi)^{\frac{N-2s}{2s}}}S^\frac{N}{2s}.
\end{eqnarray}

Now, by Lemma \ref{relation} and Remark \ref{OBS10}, we get
\begin{eqnarray*}
S^{\frac{N}{2s}}&=& \dfrac{(B^\alpha C^\beta +\xi_1 B^{2^*_s}+\xi_2 C^{2^*_s})^{\frac{N-2s}{2s}}}{(B^2+C^2)^{\frac{N}{2s}}} S_s^{\frac{N}{2s}} \\ &\leq& (1+\xi)^{\frac{N-2s}{2s}}
\dfrac{\Big[( B^{2}+ C^{2})^{\frac{2^*_s}{2}}\Big]^{\frac{N-2s}{2s}}}{(B^2+C^2)^{\frac{N}{2s}}} S_s^{\frac{N}{2s}} =
(1+\xi)^{\frac{N-2s}{2s}}
S_s^{\frac{N}{2s}},
\end{eqnarray*}
and consequently by the estimate
\eqref{SS}, we conclude that

\begin{eqnarray}
\nonumber \I_{\lambda,s}(W+s\vec{e}_\epsilon)<\dfrac{s}{N}S_s^\frac{N}{2s}.
\end{eqnarray}

Second case: Let $s>s_0.$ As in the proof of Lemma \ref{functional}, from \eqref{des}, Lemma \ref{prop2} and Lemma \ref{lemm}, we derive

\begin{eqnarray*}
\I_{\lambda,s}(W+s\vec{e}_\epsilon) &\leq& \frac{1}{2} s^2 \left(||\vec{e}_\epsilon||^2 _Y - \mu_1 ||\vec{e}_\epsilon ||^2 _{L^2 \times L^2}\right)- \int_{\Omega}F (w+s\vec{e}+U_T) dx.
\end{eqnarray*}
On the other hand,
\begin{eqnarray*}
\nonumber && F (w+s\vec{e}+U_T)  = \frac{1}{2^*_s} \Big[(sB)^{\alpha}
\Big(e_{\epsilon}+\frac{w_1+u_{r,t}}{sB}\Big)^{\alpha}_+ (sC)^{\beta}
\Big(e_{\epsilon}+\frac{w_2+v_{r,t}}{sC}\Big)^{\beta}_+ \\
&&+\xi_1 (sB)^{2^*_s}\Big(e_{\epsilon}+\frac{w_1+u_{r,t}}{sB}\Big)^{2^*_s}_+ +\xi_2 (sC)^{2^*_s}
\Big(e_{\epsilon}+\frac{w_2+v_{r,t}}{sC}\Big)^{2^*_s}_+ \Big].
\end{eqnarray*}
Using previous arguments, we get
\begin{eqnarray*}
\I_{\lambda,s}(W+s\vec{e}_\epsilon) \leq \Phi_\epsilon (s),
\end{eqnarray*}
where
\begin{eqnarray*}
\Phi_\epsilon (s) &:=& \frac{1}{2} s^2 \left(||\vec{e}_\epsilon||^2 _Y - \mu_1 ||\vec{e}_\epsilon ||^2 _{L^2 \times L^2}\right)\\ &-& \frac{s^{2^*_s}}{2^*_s} (B^\alpha C^\beta + \xi_1 B^{2^*_s}+\xi_2 C^{2^*_s}) \|e_\epsilon\|^{2^*_s}_{L^{2^*_s}(\Omega_\epsilon)}+ K s^{2^*_s}O(\epsilon^q)
\end{eqnarray*}
with $q:=\min \left\{ \frac{N-2s}{2} \; , \;p \right\} > 2s$ (because $N>6s$ and $p>2$).

Applying Lemma \ref{L3} to the function $\Phi_\epsilon$, we have by Lemmas \ref{prop2}, \ref{lemm2} and \ref{lemm} and by the choice of $B$ and $C$,

\begin{eqnarray*}
\Phi_\epsilon (s) &\leq &
\Phi_s (s_\epsilon)\\
&\leq &\frac{s}{N} \frac{\left[(B^2 + C^2)S ^{\frac{N}{2s}} + O (\epsilon^{N-2s}) - \mu_1 C e^{2s} + O (\epsilon^{N-2s}) \right]^\frac{N}{2s}}{\left[(B^\alpha C^\beta +\xi_1 B^{2^*_s} +\xi_2 C^{2^*_s}) S^{\frac{N}{2s}}+ O (\epsilon^N) + O (\epsilon^{N-2s})\right]^{\frac{N-2s}{2s}}} + O (\epsilon^q)  \\
& \leq &
\frac{s}{N} \left[ \frac{(B^2 + C^2)}{(B^\alpha C^\beta +\xi_1 B^{2^*_s} +\xi_2 C^{2^*_s})^{\frac{N-2s}{N}} }S \right]^{\frac{N}{2s}}  -\frac{\mu_1}{N}O(\epsilon^{2s})+ O (\epsilon^q).
\end{eqnarray*}
Since $q>2s$, taking $\epsilon>0$ sufficiently small, we obtain
\begin{eqnarray*}
\I_{\lambda,s}(W+s\vec{e}_\epsilon)
\leq \frac{s}{N} S_s^{\frac{N}{2s}}.
\end{eqnarray*} $\hfill \rule{2mm}{2mm}$

\subsection{The Palais-Smale condition for the functional $\I_{\lambda,s}$}

In this subsection we discuss a compactness property for the functional $\I_{\lambda,s},$ given by the Palais-Smale condition.

\begin{lemma}\label{lemalim}
If $k\geq 0$ and $\lambda_{k,s} < \mu_1 \leq \mu_2 < \lambda_{k+1,s}$. Then every $(PS)_c$ sequence of $\I_{\lambda,s}$ is bounded.
\end{lemma}
\textbf{Proof}
The Fr\'echet derivative of the functional $\I_{\lambda,s}$ is given by
\begin{align*}
  \I_{\lambda,s}'(u,v)(\phi, \psi) 
=&\langle (u,v),(\varphi,\psi)\rangle_Y  
-\int_{\Omega}(A(u,v),(\phi,\psi))_{\R^2} dx
-\int_{\Omega} (\nabla F(u+u_T, v+v_T),(\phi, \psi))_{\R^2} dx,
\end{align*}
for every $(u,v), (\phi,\psi)\in Y(\Omega).$

Let $(U_n) \subset Y(\Omega)$ be a $(PS)_c$-sequence, i.e.
satisfying $ \I_{\lambda,s} (U_n) = c+o(1)$ and $ \langle \I_{\lambda,s}'(U_n),\Psi \rangle=o(1) \|\Psi\|_Y, \; \forall \, \Psi=(\psi,\xi) \in Y(\Omega). $

Therefore
\begin{align}\label{1}
\I_{\lambda,s}(U_n)- \frac{1}{2}\I_{\lambda,s}'(U_n)U_n &=\nonumber
\frac{1}{ 2} \int_{\Omega}  (\nabla F(U_n + U_T), U_n)_{\R^2}dx -\int_{\Omega}  F(U_n + U_T)dx\\
&\leq  c+o(1)+o(1)\|U_n\|_{Y}.
\end{align}

From (\ref{1}), we get
\begin{align}
\label{3} \nonumber
&\frac{1}{2}\int_{\Omega}  (\nabla F(U_n + U_T), U_n)_{\R^2}dx-\int_{\Omega}  F(U_n + U_T)dx \\ \nonumber &=\frac{1}{2}\int_{\Omega}\left( \frac{\alpha}{\alpha+\beta}(u_n+u_{r,t})_{+}^{\alpha-1}(v_n+v_{r,t})_{+}^{\beta}u_n\right.+\xi_1 (u_n+u_{r,t})_{+}^{\alpha+\beta-1}u_n\nonumber\\
&+ \frac{\beta}{\alpha+\beta}(u_n+u_{r,t})_{+}^{\alpha}(v_n+v_{r,t})_{+}^{\beta-1}v_n\nonumber \left.+\xi_2 (v_n+v_{r,t})_{+}^{\alpha+\beta-1}v_n\right)dx
\nonumber\\
&- \frac{1}{\alpha+\beta}\int_{\Omega}\left((u_n+u_{r,t})_{+}^{\alpha}(v_n+v_{r,t})_{+}^{\beta} +\xi_1
(u_n+u_{r,t})_{+}^{\alpha+\beta} + \xi_2 (v_n+v_{r,t})_{+}^{\alpha+\beta}\right)dx\\ \nonumber
&\leq  c+o(1)+o(1)\|U_n\|_{Y}.
\end{align}

Now note that
\begin{align}
\label{03} \nonumber
&\int_{\Omega}\left( (u_n+u_{r,t})_{+}^{\alpha-1}(v_n+v_{r,t})_{+}^{\beta}u_n\right) dx
\nonumber\\
&=\int_{\Omega}\left((u_n+u_{r,t})_{+}^{\alpha-1}(u_n+u_{r,t})_{+}(v_n+v_{r,t})_{+}^{\beta} \right)dx\\ \nonumber
&-\int_{\Omega}\left((u_n+u_{r,t})_{+}^{\alpha-1}(v_n+v_{r,t})_{+}^{\beta}u_{r,t} \right)dx
\end{align} and
\begin{align}
\label{04}
&\int_{\Omega}\left( (u_n+u_{r,t})_{+}^{\alpha+\beta-1} u_n\right) dx\\ \nonumber
&=\int_{\Omega}(u_n+u_{r,t})_{+}^{\alpha+\beta}dx-\int_{\Omega}(u_n+u_{r,t})_{+}^{\alpha+\beta-1}u_{r,t} dx.
\end{align}

Substituting (\ref{03}), (\ref{04}) and expressions similar to these in (\ref{3}), yields that

\begin{equation}
\left.
\begin{array}{rcl}\label{4}
\displaystyle \int_{\Omega}(u_n+u_{r,t})_{+}^{\alpha}(v_n+v_{r,t})_{+}^{\beta}dx,
\\\displaystyle \int_{\Omega}(u_n+u_{r,t})_{+}^{\alpha+\beta}dx,
\\ \displaystyle \int_{\Omega}(v_n+v_{r,t})_{+}^{\alpha+\beta}dx
\end{array}
\right\}\leq c+o(1)+o(1)\|U_n\|_Y.
\end{equation}
Now, using  (\ref{controllo}) and that $\Psi=U_n^+ = (u_n^+ , v_n^+) \in  E_k^+$, we obtain

\begin{eqnarray}\label{quatro} \nonumber
\Big(1- \frac{\mu_2}{\lambda_{k+1,s}}\Big)\|U_n^+ \|_Y^2 &\leq&  \|U_n^+\|^2_Y-\int_{\Omega}(A U_n^+,U_n^+)_{\R^2} dx \\
&=& \int_{\Omega}  (\nabla F(U_n + U_T), U_n^+)_{\R^2}dx - \langle\I_{\lambda,s}'(U_n),(U_n^+)\rangle\\
&\leq& \int_{\Omega} F_u(U_n + U_T) |u_n^+| dx + \int_{\Omega}  F_v(U_n + U_T) |v_n^+| dx + C\|U_n^+ \|_Y.\nonumber
\end{eqnarray}
Hence,  by Remark \ref{obs10} $(iii),$ there exists a constant $K>0$ such that
\begin{align*}
F_{u}(U) &\leq K\left((u)_{+}^{\alpha+\beta-1}+(v)_{+}^{\alpha+\beta-1}\right)\\
F_{v}(U) &\leq K\left((u)_{+}^{\alpha+\beta-1}+(v)_{+}^{\alpha+\beta-1}\right).
\end{align*}

Then
\begin{align*} \nonumber
&\int_{\Omega} F_u(U_n + U_T) |u_n^+| dx + \int_{\Omega}  F_v(U_n + U_T) |v_n^+| dx  \\ \nonumber
\leq&  K\int_{\Omega}\left( (u_n + u_T)_+^{\alpha+\beta-1} +
(v_n + v_T)_+^{\alpha+\beta-1}\right)|u_n^+|dx \\ \nonumber
+& K\int_{\Omega}\left( (u_n + u_T)_+^{\alpha+\beta-1} +
(v_n + v_T)_+^{\alpha+\beta-1}\right)|v_n^+|dx.
\end{align*}
and using H\"{o}lder's inequality with $p=\frac{\alpha+\beta}{\alpha+\beta-1}$ and $q=\alpha+\beta$, Young's inequality, we deduce that
\begin{align*} \nonumber
&\int_{\Omega} F_u(U_n + U_T) |u_n^+| dx + \int_{\Omega}  F_v(U_n + U_T) |v_n^+| dx  \\ \nonumber
\leq&  K \Big\{ \epsilon \|u_n^+\|^{2}_{L^{2^*_s}}  + C_{\epsilon} \Big[ \| (u_n + u_T)_+ \|^{2(2^*_s-1)}_{L^{2^*_s}} +
\| (v_n + v_T)_+ \|^{2(2^*_s-1)}_{L^{2^*_s}} \Big] \Big\}\\ \nonumber
+& K \Big\{ \epsilon  \|v_n^+\|^{2}_{L^{2^*_s}}  + C_{\epsilon} \Big[ \| (u_n + u_T)_+ \|^{2(2^*_s-1)}_{L^{2^*_s}} +
\| (v_n + v_T)_+ \|^{2(2^*_s-1)}_{L^{2^*_s}} \Big] \Big\}.
\end{align*}

Using (\ref{4}), in view of the embedding $X(\Omega) \hookrightarrow L^r (\Omega), \; \forall \; r \leq 2^*_s$, we get
\begin{align*}
&\int_{\Omega} F_u(U_n + U_T) |u_n^+| dx + \int_{\Omega}  F_v(U_n + U_T) |v_n^+| dx  \\ \nonumber
&\leq  \epsilon C_1 \| U_n^+ \|_Y^2+C_2C_{\epsilon} + 4\epsilon_n  \|U_n \|_Y^{\frac{N+2s}{N}}.
\end{align*}
By (\ref{quatro}), taking $\epsilon>0$ small enough, we conclude that
\begin{align}\label{menos}
\|U_n^+ \|_Y^2  \leq C_3 + C_4\|U_n \|_Y^{\frac{N+2s}{N}} + C_5 \|U_n^+ \|_Y.
\end{align}
Analogously, the following estimate is valid
\begin{align}\label{mais}
\|U_n^- \|_Y^2  \leq C_6+ C_7\|U_n \|_Y^{\frac{N+2s}{N}} + C_8\|U_n^+ \|_Y.
\end{align}
Using the estimates (\ref{menos}) and  (\ref{mais}), we get
\begin{align*}
\|U_n \|_Y^2  \leq C + C\|U_n \|_Y^{\frac{N+2s}{N}} + C \|U_n\|_Y.
\end{align*}
Since $\dfrac{N+2s}{N} <2,$ we conclude that $(U_n)$ is bounded in $Y(\Omega).$
$\hfill \rule{2mm}{2mm}$

\begin{lemma}\label{lema10}
If $k\geq0$ and $\lambda_{k,s} < \mu_1 \leq \mu_2 < \lambda_{k+1,s}$, then the functional $\I_{\lambda,s}$ satisfies the $(PS)$ condition at level $c$ with $c<\dfrac{s}{N}S_s^\frac{N}{2s}$.
\end{lemma}
\textbf{Proof}
Let $(U_n)\subset Y(\Omega)$ be a sequence satisfying
\begin{equation*}
\I_{\lambda,s}(U_n)\to c\quad\mbox{and}\quad \I_{\lambda,s}'(U_n)\rightarrow 0\quad \mbox{in the dual space}\; Y(\Omega)',
\end{equation*}
as $n\to\infty.$ By Lemma \ref{lemalim} we have that $(U_n)$ is bounded. Hence passing to a subsequence, we may suppose that
\begin{align}\label{sub}
U_n &\rightharpoonup U \quad \mbox{in} \;  Y(\Omega),\nonumber\\
U_n &\rightarrow U \quad \mbox{in} \; L^p (\Omega) \times L^p (\Omega), \;\mbox{for all} \; p \in [1, 2^*_s),\\
U_n &\rightarrow U \quad \mbox{a.e}\quad  \mbox{in}\; \R^N.\nonumber
\end{align}

Hence, $U$ is s weak solution to
\begin{equation}\label{sf}
\left\{
\begin{aligned}
(-\overrightarrow{\Delta})^s U &= AU + \nabla F(U+U_T) && \text{in $\Omega$,} \\
U = & \;0  && \text{in $\R^N\setminus\Omega$},
\end{aligned}
\right.
\end{equation}
that is, for any $\Psi\in Y(\Omega)$ it holds
\begin{align}\label{sfd}
\langle U,\Psi\rangle_Y
-\int_{\Omega}(AU,\Psi)_{\R^2}&=\int_{\Omega}(\nabla F(U+U_T),\Psi)_{\R^2} dx.
\end{align}
In particular, taking $\Psi=U$ in \eqref{sfd},we get
\begin{align}\label{s}
\|U\|_Y^2-\int_{\Omega}(AU,U)_{\R^2} dx=&\int_{\Omega}(\nabla F(U+U_T),U)_{\R^2} dx
\end{align}
Note that by \eqref{sf} $\left(\langle \I_{\lambda,s}'(U),U \rangle=0\right)$ and \eqref{s} we obtain
\begin{align}\label{posi}
\I_{\lambda,s}(U)= \frac{1}{ 2} \int_{\Omega}  (\nabla F(U + U_T), U)_{\R^2}dx -\int_{\Omega}  F(U + U_T)dx\geq 0.
\end{align}

By applying the Brezis-Lieb Lemma \cite{bl}, it follows that
\begin{align} \label{bl}
\| (U_n + U_T)_+ \|^{2^*_s}_{{L^{2^*_s}}\times{L^{2^*_s}}}& = \| (U_n - U)_+ \|^{2^*_s}_{{L^{2^*_s}}\times{L^{2^*_s}}}+\| (U + U_T)_+ \|^{2^*_s}_{{L^{2^*_s}}\times{L^{2^*_s}}}+o(1)\nonumber\\
\|U_n - U\|^{2}_{Y}& = \| U_n \|^{2}_{Y}-\|U \|^{2}_{Y}+o(1)
\end{align}
and by applying the Brezis-Lieb Lemma for homogeneous functions \cite{MD}, we conclude that
\begin{align}\label{blh}
\int_{\Omega}F(U_n+U_T)dx= \int_{\Omega}F(U+U_T)dx+ \int_{\Omega}F(U_n-U)dx+o(1).
\end{align}
Also, we have
\begin{align}\label{de}
\int_{\Omega}(\nabla F(U_n+U_T),U_n+U_T)_{\R^2} dx &- \int_{\Omega}(\nabla F(U+U_T),U+U_T)_{\R^2}  dx\nonumber\\&= (\alpha+\beta)\int_{\Omega} F(U_n-U) dx.
\end{align}
Then, by using $\eqref{sub}, \eqref{bl}$ and  $\eqref{de},$ we deduce
\begin{align}\label{func}
\I_{\lambda,s}(U_n)=\dfrac{1}{2}\|U_n-U\|^{2}_{Y} +\I_{\lambda,s}(U)-\int_{\Omega}F(U_n-U)dx+o(1).   \end{align}
On the other hand, by using $\eqref{sub}, \eqref{s}$ and  $\eqref{bl}$ and $\eqref{blh}$ we have
\begin{align*}
\langle&\I_{\lambda,s}'(U_n),U_n \rangle=\|U_n\|^{2}_{Y}-\int_{\Omega}(AU_n,U_n)_{\R^2} dx -\int_{\Omega}(\nabla F(U_n+U_T),U_n)_{\R^2} dx\\
=&\left[\|U_n-U\|^{2}_{Y}+\|U\|^{2}+o(1)\right]-\left[\int_{\Omega}(AU,U)_{\R^2} dx +o(1)\right]\\&-\int_{\Omega}(\nabla F(U_n+U_T),U_n+U_T)_{\R^2} dx
+\int_{\Omega}(\nabla F(U_n+U_T),U_T)_{\R^2} dx\\
=&\|U_n-U\|^{2}_{Y}+\left[\|U\|^{2}_{Y}-\int_{\Omega}(AU,U)_{\R^2} dx \right]
-\left[(\alpha+\beta)\int_{\Omega} F(U_n-U)dx \right.\\&+\left.\int_{\Omega}(\nabla F(U+U_T),U+U_T)_{\R^2}  dx\right]+\int_{\Omega}(\nabla F(U_n+U_T),U_T)_{\R^2}  dx +o(1)\\
=&\|U_n-U\|^{2}_{Y}+\left[\int_{\Omega}(\nabla F(U+U_T),U)_{\R^2}  dx \right]
-\left[(\alpha+\beta)\int_{\Omega} F(U_n-U)dx \right.\\&\left.+\int_{\Omega}(\nabla F(U+U_T),U)_{\R^2} dx +\int_{\Omega}(\nabla F(U+U_T),U_T)_{\R^2} dx \right]+\int_{\Omega}(\nabla F(U_n+U_T),U_T)_{\R^2} dx +o(1)\\
=&\|U_n-U\|^{2}_{Y}
-(\alpha+\beta)\int_{\Omega} F(U_n-U) dx+\int_{\Omega}(\nabla F(U+U_T),U_T)_{\R^2} dx \\&+\int_{\Omega}(\nabla F(U_n+U_T),U_T)_{\R^2} dx +o(1).
\end{align*}
Taking into account that $\langle\I_{\lambda,s}'(U_n),U_n \rangle \rightarrow 0$ and $\displaystyle\int_{\Omega}(\nabla F(U_n+U_T),U_T)_{\R^2} dx \rightarrow\int_{\Omega}(\nabla F(U+U_T),U_T)_{\R^2} dx $ as $n\rightarrow\infty,$ we deduce
\begin{align}\label{lim}
\|U_n-U\|^{2}_{Y}=(\alpha+\beta)\int_{\Omega} F(U_n-U)dx +o(1).
\end{align}
Let $$L:=\lim_{n\to\infty}\|U_n-U\|^{2}_{Y}\geq 0.$$
If $L=0,$ then $U_n\to U$ in $Y(\Omega)$ as $n\to\infty.$

Let $L>0.$
Then, by the definition of $S_s$,
$$S_s\leq\dfrac{\|U\|^{2}_{Y}}{\left(\displaystyle\int_{\Omega} |u|^\alpha|v|^\beta+\xi_1|u|^{\alpha+\beta}+\xi_2|v|^{\alpha+\beta}dx\right)^{\frac{2}{\alpha+\beta}}}\quad \mbox{for all}\;U=(u,v)\neq
(0,0).$$ and \eqref{lim}, we can infer
\begin{align*}
\|U_n-U\|^{2}_{Y}
&\geq S_s\left(\displaystyle\int_{\Omega} (u_n-u)_{+}^\alpha(v_n-v)_{+}^\beta+\xi_1(u_n-u)_{+}^{\alpha+\beta}+\xi_2(v_n-v)_{+}^{\alpha+\beta}dx\right)^{\frac{2}{\alpha+\beta}}\\
&=S_s\left((\alpha+\beta)\int_{\Omega} F(U_n-U)dx \right)^{\frac{2}{\alpha+\beta}}
\end{align*}
which gives
\begin{align}\label{so}
L\geq S_s L^{\frac{N-2s}{N}},\;\mbox{i.e.}\; L\geq S_s^{\frac{N}{2s}}.
\end{align}
Now, from $\eqref{posi}, \eqref{func},\eqref{lim},\eqref{so}$
we get
$$\dfrac{s}{N}S_s^{\frac{N}{2s}}\leq\left( \dfrac{2s}{N-2s}\right) \dfrac{L}{2^*_s}\leq c<\dfrac{s}{N}S_s^{\frac{N}{2s}},$$ which contradiction.
$\hfill \rule{2mm}{2mm}$

{\bf Proof of Theorem \ref{teo2}.}
In the case where $\lambda_{k,s}<\mu_1 \leq\mu_2 < \lambda_{k+1,s}$ occurs,
 the Proposition \ref{GMPG} and Lemma \ref{functional}  with $\varepsilon >0$
small enough,
ensure that the functional $\I_{\lambda,s}$ satisfies the geometric structure required by the Linking Theorem. Therefore, it follows from the Linking Theorem without the Palais-Smale condition, that there exists a sequence $(U_n) \subset Y(\Omega)$
satisfying $ \I_{\lambda,s} (U_n) \rightarrow c$ and $ \I_{\lambda,s}'(U_n) \rightarrow 0 $ in $Y(\Omega)',$
and by Lemma \ref{mini}, the critical level satisfies
$$0<c:= \inf_{\gamma \in \Gamma} \sup_{U \in Q} \I_{\lambda,s}(\gamma(U))\leq \frac{s}{N} S_s^{\frac{N}{2s}},$$
where $\Gamma:= \{ \gamma \in C^0 (Q, Y(\Omega)): \; \gamma=Id \; \mbox{on} \; \partial Q \}.$
By Lemma \ref{lemalim}, $(U_n)$ is bounded in $ Y(\Omega)$ and consequently the Lemma \ref{lema10} ensures that $U_n \rightarrow \overline{U}$ in $Y(\Omega)$.\\
If $ 0 = \lambda_{0,\mu}  <  \mu_1  \leq \mu_2 < \lambda_{1,\mu}$, to show that the functional $\I_{\lambda,s}$
satisfies the geometrical conditions of the Mountain Pass Theorem, it is enough to take the finite dimensional subspace $E^{-}=\{ (0,0) \}$ and to apply the Proposition \ref{GMPG} with $E^{+}_k=Y(\Omega)$  such that $R \| \vec{e}_{\epsilon} \|_Y > \rho$ with $R>0$ sufficient large to ensure that $\I_{\lambda,s}(R \vec{e}_{\epsilon} )<0.$ The $(PS)_c$ condition is guaranteed by making $k = 0$ in the Lemmas \ref{lemalim} and \ref{lema10}.
Thus, in both cases,
there exists a non-trivial solution $\overline{U}$ for the problem (\ref{c1.7}).
By \cite{frp} Remark 4.1, it follows that $\overline{U}_+ \neq 0$ and therefore, $U_T$ and $U_T + \overline{U}$ are distinct solutions for the problem (\ref{c1.1}). $\hfill \rule{2mm}{2mm}$

\section{The resonant case}

\subsection{Proof of Theorem \ref{teo4}}

In this subsection we discuss a compactness property for the functional $\I_{\lambda,s},$ given by the Palais-Smale condition for this case.

\begin{lemma}\label{lemreso}
If $N>6s$ and $\lambda_{k,s} = \mu_1 \leq \mu_2 < \lambda_{k+1,s}$ for $k > 1$, the functional $\I_s$ satisfies the $(PS)$ condition.
\end{lemma}
\textbf{Proof} We follow the notations of the previous proof.\\
Let $U_n \in Y(\Omega)$ such that $\I_s(U_n)\rightarrow c$ and $\I_s'(U_n) \rightarrow 0$ in the dual space $Y(\Omega)' .$
Writing  $Y(\Omega)= E^-_{k-1} \oplus E^+_k \oplus Z_k,$ consequently we have
$$U_n=U_n^- +U_n^+ +\beta_n Y_n:=W_n+\beta_n Y_n,$$ where $U_n^- \in E^-_{k-1}$, $U_n^+ \in E^+_k=(E^-_k)^\perp$
and $ Y_n \in Z_k=span \{(\varphi_{k,s},0), (0,\varphi_{k,s}) \}$ with $\|Y_n\|_Y =1$.\\
Using similar arguments as in (\ref{menos}) and (\ref{mais}), we obtain
\begin{equation}\label{eqp5}
 \|W_n \|_Y^2  \leq    C+ C\|U_n \|_Y^{\tau} +  C \|W_n\|_Y,
 \end{equation}
where $\tau=\frac{N+2s}{N}.$
We can assume $\|U_n\|_Y \ge 1$ (if $\|U_n\|_Y \le 1$, the sequence $(U_n)$ is bounded in $Y(\Omega)$). Then, since $\|U_n\|_Y \le \|W_n \|_Y + |\beta_n|$, from (\ref{eqp5}), we have
\begin{equation}
\| W_n \|_Y^2\le C_1 ( \| W_n \|_Y+ | \beta_n |)^{\tau}+ C\|W_n\|_Y. \label{pstar}
\end{equation}
If $\beta_n$ is bounded, since $\tau <2,$ by (\ref{pstar}) we conclude that $(U_n)$ is bounded in $Y(\Omega)$.
Otherwise, we may assume $\beta_n \to +\infty,$ therefore, from (\ref{pstar}), it follows that
\begin{equation}\nonumber
 \left\| \frac{W_n}{\beta_n} \right\|_Y^2\le C_1  \left\{ \frac{( \| W_n \|_Y+| \beta_n|)^{\tau/2}}{| \beta_n|} \right\}^2+
 C \frac{1}{\beta_n}  \left\| \frac{W_n}{\beta_n} \right\|_Y
\end{equation}
\begin{equation}\label{XX}
 \leq  C_1  \left\{   \frac{1}{|\beta_n|^{1-\tau/2}} \left\| \frac{W_n}{\beta_n} \right\|_Y^{\tau/2} + \frac{1}{| \beta_n|^{1-\tau/2}} \right\}^2+
 C \frac{1}{\beta_n}  \left\| \frac{W_n}{\beta_n} \right\|_Y.
\end{equation}
Using again the fact that $\tau/2 <1,$ the above estimate yields that
\begin{equation}\nonumber
  \left\| \frac{W_n}{\beta_n} \right\|_Y^2\leq  C_2  \left\| \frac{W_n}{\beta_n} \right\|_Y^{\tau} +
 C_3  \left\| \frac{W_n}{\beta_n} \right\|_Y + C_4
\end{equation}
and consequently the sequence $\left\{\dfrac{W_n}{\beta_n} \right\}$ is bounded in $Y(\Omega)$ and by \eqref{XX},
$ \Big\| \dfrac{W_n}{\beta_n} \Big\|_Y \to 0.$
Therefore, possibly up to a subsequence,
$W_n/\beta_n\to 0$ a.e. in $\Omega$ and strongly in $L^q(\Omega) \times L^q(\Omega) $, $1\le q<2^*_s$;
$Y_n\to Y_0 \in Z_k$ a.e. in $\Omega$ and strongly in $Y(\Omega)$ and $L^q(\Omega) \times L^q(\Omega) $, $1\le q<2^*_s$.\\
Now, taking $\beta_n Y_n \in Z_k$ as test function, we get
$$ \I_s'(U_n)Y_n = \beta_n \left(\| Y_n\|_Y^2 -  \int_{\Omega} (AY_n,Y_n)_{\R^2} dx\right)  - \int_{\Omega} (\nabla F(U_n + U_T),Y_n)_{\R^2} dx.$$
Since $(U_n)$ is a $(PS)$-sequence and $\frac{1}{(\beta_n)^{\frac{4s}{N-2s}}} \left(\| Y_n\|_Y^2 -  \displaystyle\int_{\Omega} (AY_n,Y_n)_{\R^2} dx\right)\rightarrow 0,$
as $n\rightarrow \infty,$ we obtain that
$$o(1)=  \frac{1}{(\beta_n)^{\frac{N+2s}{N-2s}}} \I_s'(U_n) (Y_n) =  -  \frac{1}{(\beta_n)^{\frac{N+2s}{N-2s}}}\int_{\Omega} (\nabla F(U_n + U_T),Y_n)_{\R^2} dx.$$

Now, from Remark \ref{obs10} $(ii)$,
\begin{equation} \displaystyle\label{nabla} \int_{\Omega} (\nabla F\Big(\frac{U_n + U_T}{\beta_n}\Big),Y_n)_{\R^2} dx= \frac{1}{(\beta_n)^{\frac{N+2s}{N-2s}}}
 \int_{\Omega} (\nabla F(U_n + U_T),Y_n)_{\R^2} dx \rightarrow 0.
 \end{equation}

 On the other hand, since $U_n=W_n+\beta_n Y_n,$ we have that $\dfrac{U_n}{\beta_n} \rightarrow Y_0$ in $L^q(\Omega) \times L^q(\Omega)$ for all $1 \leq q < 2^*_s$
and a.e in $\Omega.$ So, by the Dominated Convergence Theorem and by (\ref{nabla}), it follows that
$$\int_{\Omega} (\nabla F\Big(\frac{U_n + U_T}{\beta_n}\Big),Y_n)_{\R^2} dx \rightarrow \int_{\Omega} (\nabla F(Y_0),Y_0)_{\R^2} dx=0$$
and from Remark \ref{obs10} $(i)$, we concluded that $\displaystyle \int_{\Omega} F(Y_0)dx=0.$\\
Finally, using the notation $Y_0= (y_1^0, y_2^0)$, it follows that $(y_1^0)_+ = 0= (y_2^0)_+$,
contradicting $\| Y_0 \|_Y = 1$ and $Y_0 \in Z_k$ with $k>1$, which ensures that at least one of the functions is not null and changes sign.
Thus  $(U_n)$ is bounded and using the fact that $N>6s,$ as in the proof of Lemmas \ref{mini} and \ref{lema10}, we have that $(U_n)$ admits a convergent subsequence.
$\hfill \rule{2mm}{2mm}$

\subsection{Geometry in resonant case}

In this subsection, we demonstrate that the functional $\I_{\lambda,s}$ satisfies the geometric structure required by the Linking Theorem in resonant case, that is, we obtain the following result.
\begin{proposition}\label{rest}
Suppose $\Omega$ is a smooth bounded domain of $\R^N$, $\alpha+\beta=2^{*}_{s}$ and $\lambda_{k,s}= \mu_1 \leq \mu_2< \lambda_{k+1,s}$ for some $k>1.$ Then
\item[i)] there exist $\sigma, \rho>0$ such that $\I_s(U)\geq \sigma$ for all $U \in E^+_k$ with $\|U\|_{Y}=\rho,$
\item[ii)] there exists $E \in  E_k^{+}$ and $R>0$ such that $R \|E\|_Y > \rho$ and $\I_s(U) \leq 0,$ for all $U \in \partial Q,$ where
$Q=(\overline{B}_{R}\cap E^-_{k})\oplus[0,R]E.$
\end{proposition}
\textbf{Proof}
\textbf{i)} Let $U=(u,v) \in E^+_k$, using the fact that $u_T, v_T < 0$, estimate $|u|^{\alpha} |v|^{\beta} \leq |u|^{\alpha+\beta} + |v|^{\alpha+\beta} $
and the fractional imbedding $X \hookrightarrow L^{\alpha+\beta}$, by (\ref{controllo}), we have
\begin{eqnarray*}\I_s(U)&\geq& \frac{1}{2} \| U\|_Y^2 - \frac{\mu_2}{2}  \| U\|_{(L^2)^2}^2 - C \int_{\Omega} ( | u|^{\alpha+\beta} + | v |^{\alpha+\beta} )dx\\ &\geq&
\frac{1}{2} \Big( 1 - \frac{\mu_2}{\lambda_{k+1,s}} \Big)\|U\|^{2}_{Y}
- C\|U\|^{\alpha+\beta}_{Y},\end{eqnarray*}
where $C>0$ is a constant.
Since $\mu_2 < \lambda_{k+1,s}$ and $\alpha+\beta>2,$ for $ \| U\|_Y =\rho$ small enough, we get $\I_s(U)\geq \sigma.$

\textbf{ii)} Now consider the following decomposition
$\begin{array}{l}
\partial Q=\Gamma_1 \cup \Gamma_2 \cup \Gamma_3,
\end{array}
$
where
$
\begin{array}{l}
 \Gamma_1=\{U\in Y(\Omega);\; U=U_1+ r E,
 \; \mbox{with}\; U_1 \in E^{-}_k,\; \|U_1\|_Y=R, \; 0 \leq r \leq R\}, \\
 \Gamma_2=\{U\in Y(\Omega); \; U=U_1+ R E,\
 \; \mbox{with}\; U_1 \in E^{-}_k,\; \|U_1\|_Y\leq R \},\\
 \Gamma_3= \overline{B_{R}(0)}\cap E^-_k.
\end{array}
$

Let us show that on each set $\Gamma_i$ we have $\I_s \mid_{\Gamma_i}\leq 0, \; i=1,2,3.$

Choose $E$ as follows:\\
Fixed $R_0 > \rho$, take $E=(e_1, e_2) \in E_k^+ = (E_k^-)^{\perp}$ (with $e_i \geq 0, \; i=1,2$) 
satisfying\\
(I) $\| E \|^2_Y <\Big( \dfrac{\mu_1}{\lambda_{k-1,s}} -1 \Big)\delta^2,$ where $\delta>0$ is a constant
to be obtained forward.\\
(II) $e_1 \geq 2 \Big( K + \frac{\| u_T \|_{C^0}}{R_0}  \Big)$ and $e_2 \geq 2 \Big( K + \frac{\| v_T \|_{C^0}}{R_0}  \Big)$ a.e. in some 
$\mathcal{C} \subset \Omega$ with $|\mathcal{C}|>0,$
where $K>0$ satisfies $\| V\|_{(C^0)^2} \leq K \| V\|_Y,$ for all $V \in E_k^-.$\\ Note that this choice is possible because $(E_k^-)^{\perp}$ has unbounded functions; 
$E_k^-$ has finite dimension and $K=  \sup_{\scriptsize{ \begin{array}{c}
\| V\|_Y=1 \\
V \in E_k^- \\
\end{array}} } \| V\|_{(C^0)^2}.$

\textbf{ Estimates on $\Gamma_1$:}
 For $U=U_1+rE \in \Gamma_1$, we consider $U_1=R\widehat{U}_1 \in E_{k}^-$ with $\|\widehat{U}_1\|_E=1$ and we set $\widehat{U}_1=c_1Y+c_2E_k,
\mbox{ where }E_k\in Z_k=\mbox{ span}\{(\varphi_{k,s},0),(0,\varphi_{k,s})\}$ and $Y\in E_{k-1}^{-}$ with $\| Y \|_Y=1$.
Then, 
\begin{equation}
\begin{array}{lll}
 \I_s(U) &\leq& \dfrac{1}{2} \| U_1 \|^2_Y + \dfrac{r^2}{2} \|E\|^2_Y -\dfrac{\mu_1}{2} \| U_1\|^2_{{(L^2)}^2} - \displaystyle \int_{\Omega} F(U+U_T) dx \\ \nonumber
 &\leq& \dfrac{R^2}{2} \| \widehat{U}_1 \|^2_Y + \dfrac{R^2}{2} \|E\|^2_Y -\dfrac{\mu_1 R^2}{2} \| \widehat{U}_1\|^2_{{(L^2)}^2}  -\displaystyle  \int_{\Omega} F(U+U_T) dx\\ \nonumber
 &=& \dfrac{R^2}{2} \| c_1Y+c_2E_k \|^2_Y + \dfrac{R^2}{2} \|E\|^2_Y -\dfrac{\mu_1 R^2}{2} \| c_1Y+c_2E_k\|^2_{{(L^2)}^2}\\
 &-& \displaystyle \int_{\Omega} F(U+U_T) dx\\ \nonumber
 &=& \dfrac{R^2}{2} c_1^2 ( \| Y \|_Y^2 - \mu_1 \|Y\|^2_{{(L^2)}^2}) +  \dfrac{R^2}{2} c_2^2 ( \|E_k\|^2_Y -\mu_1 \|E_k\|^2_{{(L^2)}^2})\\
 &+& \dfrac{R^2}{2} \|E\|^2_Y - \displaystyle \int_{\Omega} F(U+U_T) dx. \nonumber
\end{array}
\end{equation} 
Consequently
\begin{equation}\label{xxxzz}
\begin{array}{lll}
\I_s(U) \leq \dfrac{R^2}{2} c_1^2 \Big( 1 - \dfrac{\mu_1}{\lambda_{k-1,s}} \Big)\| Y \|_Y^2 +  \dfrac{R^2}{2} \|E\|^2_Y  -\displaystyle  \int_{\Omega} F(U+U_T) dx. 
\end{array}
\end{equation}

Now using the notation $\widehat{U}_1 = (\widehat{u}_1,\widehat{v}_1)=(c_1 y_1 + c_2 e_1^k, c_1 y_2 + c_2 e_2^k)$, where $Y=(y_1, y_2) \in E_{k-1}^- \cap B_1$ 
and  $E_k=(e_1^k, e_2^k) \in Z_k \cap B_1$, we will prove that there exist $\delta > 0$ and $\eta > 0$ such that
$$ \max_{i=1,2} \Big\{ \max_{\overline{\Omega}}
\{ c_1 y_i + c_2 e_i^k; \; \; | c_1 |\leq\delta \} \Big\} \ge \eta > 0.$$
Indeed, by contradiction, assume that there exist sequences $(c_1^n), (c_2^n) \subset \R$ and $Y_n= (y_1^n, y_2^n) \subset Y(\Omega)$ with $\| Y_n \|_Y=1$ such that
$c_1^n \rightarrow 0, \; |c_2^n|=\sqrt{1-(c_1^n)^2}\rightarrow 1$ and $$\max_{i=1,2} \Big\{ \max_{\overline{\Omega}}
\{ c_1^n y_i^n + c_2^n e_i^k \} \Big\} \rightarrow 0,\; \mbox{as} \; n \rightarrow \infty.$$
Therefore, $c_1^n y_i^n \rightarrow 0$ and $c_2^n e_i^k \rightarrow e_i^k$ and consequently
$$\max_{i=1,2} \Big\{ \max_{\overline{\Omega}} e_i^k (x) \Big\}=0.$$
Hence, we conclude that $e_1^k \leq 0 $ and $e_2^k \leq 0$ in $\Omega$, which is a contradiction, because $k>1,$ $E_k=(e_1^k, e_2^k) \in Z_k $ and $\| E_k \|_Y =1$ imply that at least one of the coordinate functions must change sign.\\
So, we conclude that there exist $ \delta > 0$, $\eta > 0$ such that
$$\max \Big\{  \max_{\overline{\Omega}} \widehat{u}_1; \; \max_{\overline{\Omega}} \widehat{v}_1 \;: \, | c_1 |\leq\delta  \Big\} \ge \eta > 0 \ , \
\forall \; { \widehat{U}_1}= c_1 Y + c_2 E_k \in E_{k}^- \; \mbox{with} \; \|\widehat{U}_1\|_Y = 1.$$

Denoting $\Omega_{+}=\Big\{ x\in \overline{\Omega}:  (\widehat{u}_1) (x)  \ge \eta/2 \; \mbox{and} \; (\widehat{v}_1)(x) \ge \eta/2 \Big\}$. By equicontinuity of the functions $\widehat{U}_1,$ we
have that $|\Omega_{+}| \ge \nu > 0$, $\forall {\widehat{U}_1} \in E_{k}^- \cap B_1$ and $|c_1|\le \delta.$\\
Moreover
\[
\frac{u_T(x)}{R}\geq -\frac{\|u_T\|_{C^0}}{R}>-\frac{\eta}{4} \; \mbox{and} \; \frac{v_T(x)}{R}\geq -\frac{\|v_T\|_{C^0}}{R}>-\frac{\eta}{4}, \; \forall \; R\geq R_0 \mbox{ sufficiently
large}.
\]
Then, since $e_1, e_2 \geq 0$ in $\Omega,$
\begin{eqnarray}
\nonumber \int_{\Omega} F(U+U_T) dx &\geq& \frac{\xi_1}{\alpha+\beta}  R^{\alpha+\beta}\int_{\Omega} \Big(\widehat{u}_1 + \frac{u_T}{R}\Big)_+^{\alpha+\beta} dx\\ \nonumber
&+& \frac{\xi_2}{\alpha+\beta}  R^{\alpha+\beta}\int_{\Omega} \Big(\widehat{v}_1 + \frac{v_T}{R}\Big)_+^{\alpha+\beta} dx\\ \nonumber
&\geq& C R^{\alpha+\beta} \Big[\int_{\Omega_+} \Big(\widehat{u}_1 - \frac{\eta}{4}\Big)_+^{\alpha+\beta} dx + 
 \int_{\Omega_+} \Big(\widehat{v}_1 - \frac{\eta}{4}\Big)_+^{\alpha+\beta} dx\Big]  \\ \nonumber
& \geq& C R^{\alpha+\beta} \Big[\int_{\Omega_+} \Big( \frac{\eta}{4}\Big)^{\alpha+\beta} dx + 
 \int_{\Omega_+} \Big( \frac{\eta}{4}\Big)^{\alpha+\beta} dx\Big]  \\ \nonumber
& \geq& C R^{\alpha+\beta} \Big( \frac{\eta}{4}\Big)^{\alpha+\beta} |\Omega_+|= \widetilde{C} R^{\alpha+\beta}, \nonumber
\end{eqnarray}
for all $ R$ sufficiently large.
Thus, from (\ref{xxxzz}) we can conclude that there exists $R_1 > 0$ such that
\[
\I_s(U) \leq \dfrac{R^2}{2} \delta^2 \Big( 1 - \frac{\mu_1}{\lambda_{k-1,s}} \Big) +  \dfrac{R^2}{2} \|E\|^2_Y - \widetilde{C} R^{\alpha+\beta}<0,
\]
for all $R\ge R_1$. 

On the other hand, if $ | c_1 |\geq \delta>0,$ by the choose of $E$, we get 
\begin{eqnarray}
\nonumber \I_s(U) &\leq& -\dfrac{R^2}{2} c_1^2 \Big(  \dfrac{\mu_1}{\lambda_{k-1,s}} -1\Big) +  \dfrac{R^2}{2} \|E\|^2_Y \\ \nonumber
 &\leq& -\dfrac{R^2}{2} \Big[ \delta^2  \Big(  \dfrac{\mu_1}{\lambda_{k-1,s}} -1\Big) -  \|E\|^2_Y \Big]<0. \nonumber
\end{eqnarray} 
\textbf{ Estimates on $\Gamma_2$:}
For $U=U_1+RE \in \Gamma_2$, we have
\begin{eqnarray}
\nonumber \I_s(U_1+RE) \leq  \displaystyle\frac{1}{2} \|U_1\|_Y^2 \Big(1- \displaystyle
\frac{\mu_1}{ \lambda_{k,s}} \Big) + \frac{R^2}{2} \|E\|^2_Y - \displaystyle  \int_{\Omega} F(U_1+RE+U_T) dx. \nonumber 
\end{eqnarray}
Since $\lambda_{k,s}=\mu_1$,
\begin{eqnarray}\label{ff}
\I_s(U_1+RE) \leq  \frac{R^2}{2} \|E\|^2_Y - \displaystyle  \int_{\Omega} F(U_1+RE+U_T) dx. 
\end{eqnarray}
Now, to estimate the last integral, note that, if $U_1=(u_1,u_2),$
\begin{align*}  & \displaystyle \int_{\Omega} F(U_1+RE+U_T) dx \nonumber \\
&\geq  \dfrac{1}{\alpha+\beta} \Big[ \xi_1 R^{\alpha+\beta} \displaystyle \int_{\Omega} \Big(e_1+\frac{u_1+u_T}{R} \Big)_+^{\alpha+\beta} dx +
\xi_2 R^{\alpha+\beta} \displaystyle \int_{\Omega} \Big(e_2+\frac{u_2+v_T}{R} \Big)_+^{\alpha+\beta} dx\Big]\nonumber
\end{align*} for $ R \geq R_0$, and by (II) each integral on the right can be estimated as follows
\begin{eqnarray}
\nonumber \int_{\Omega} \Big(e_i +\frac{u_i+w_T}{R} \Big)_+^{\alpha+\beta} dx &\geq&
 \int_{\Omega} \Big(e_i-\frac{\|u_i\|_{C^0} +\|w_T\|_{C^0}}{R} \Big)_+^{\alpha+\beta} dx\\ \nonumber
 &\geq&  \int_{\Omega} \Big(e_i-\Big(K+\frac{\|w_T\|_{C^0}}{R_0} \Big) \Big)_+^{\alpha+\beta} dx\\  \nonumber
 &\geq & \int_{\mathcal{C}} \Big(K+\frac{\|w_T\|_{C^0}}{R_0} \Big)^{\alpha+\beta} dx = \Big(K+\frac{\|w_T\|_{C^0}}{R_0} \Big)^{\alpha+\beta} |\mathcal{C}|, 
\end{eqnarray}
for $i=1,2 $ and $w_T \in \{u_T \; , \; v_T \}$.\\
Therefore, by (\ref{ff}) and by above estimates,
\begin{eqnarray}
\nonumber \I_s(U_1+RE) &\leq & \frac{R^2}{2} \|E\|^2_Y - \displaystyle c_1 R^{\alpha+\beta} \int_{\Omega}  \Big(e_1+\frac{u_1+u_T}{R} \Big)_+^{\alpha+\beta} dx\\ &-&
c_2 R^{\alpha+\beta} \int_{\Omega}\Big(e_2+\frac{u_2+v_T}{R} \Big)_+^{\alpha+\beta} dx \leq  \frac{R^2}{2} \|E\|^2_Y - C R^{\alpha+\beta}.\nonumber
\end{eqnarray}
Since $\alpha+\beta>2,$ for $R\geq R_0$ we have $\I_s (U) <0,$ for all $U \in \Gamma_2.$

\textbf{Estimates on $\Gamma_3$:}
For $U\in \Gamma_3$, it follows the estimate
\begin{eqnarray}\nonumber
 \I_s(U) \leq \displaystyle
   \frac{1}{2}
 \| U \|^2_Y
 - \displaystyle \frac{\mu_1}{2} \| U \|^2_{(L^2)^2} \leq \displaystyle \frac{1}{2} \Big( 1- \frac{\mu_1}{\lambda_{k,s}} \Big) \| U \|^2_Y = 0.\nonumber
  \end{eqnarray}
  
Therefore, for all $R \geq R_0>0$, follows that $\I_s(U)\leq 0$ for all $U \in \partial Q$, concluding the desired result. $\hfill \rule{2mm}{2mm}$

{\bf Proof of Theorem \ref{teo4}.}

With the previous results, we conclude the proof of Theorem \ref{teo4} by a direct application of the Linking Theorem and arguing as in proof of Theorem \ref{teo2} to obtain two distinct solutions for the problem (\ref{c1.1}). $\hfill \rule{2mm}{2mm}$

\medskip

\end{document}